\documentclass[12pt,leqno,draft]{article}
\usepackage{amsfonts}
\usepackage{amsmath, amsthm, amsfonts, amssymb, color}
\usepackage{mathrsfs}
\usepackage{color}
\theoremstyle{definition}

\newcommand{\scr}[1]{\mathscr #1}
\definecolor{wco}{rgb}{0.5,0.2,0.3}

\numberwithin{equation}{section} \theoremstyle{remark}

\newcommand{\ua}{\uparrow}

\title{{\bf
Entropy-Cost Inequalities for McKean-Vlasov SDEs with Singular Interactions}\footnote{Supported in part by the National Key R\&D Program of China (2022YFA1006000), NNSFC(12301180, 12271398),  Deutsche Forschungsgemeinschaft (DFG, German Research Foundation) (Project-ID 317210226, SFB 1283). Panpan Ren is supported by Research Centre for Nonlinear Analysis at Hong Kong PolyU.} }
\author{
{\bf   Xing Huang$^{a)}$, Panpan Ren$^{b)}$,  Feng-Yu Wang$^{a)}$  }\\
\footnotesize{ a) Center for Applied Mathematics and KL-AAGDM, Tianjin
University, Tianjin 300072, China}\\
\footnotesize{ b) Department of Mathematics, City University of  Hong Kong, Tat Chee Av., Hong Kong,  China }\\
\footnotesize{  xinghuang@tju.edu.cn, panparen@cityu.edu.hk, wangfy@tju.edu.cn}\\
}
\begin{document}
\allowdisplaybreaks
\def\R{\mathbb R}  \def\ff{\frac} \def\ss{\sqrt} \def\B{\mathbf
B} \def\W{\mathbb W}
\def\N{\mathbb N} \def\kk{\kappa} \def\m{{\bf m}}
\def\ee{\varepsilon}\def\ddd{D^*}
\def\dd{\delta} \def\DD{\Delta} \def\vv{\varepsilon} \def\rr{\rho}
\def\<{\langle} \def\>{\rangle} \def\GG{\Gamma} \def\gg{\gamma}
  \def\nn{\nabla} \def\pp{\partial} \def\E{\mathbb E}
\def\d{\text{\rm{d}}} \def\bb{\beta} \def\aa{\alpha} \def\D{\scr D}
  \def\si{\sigma} \def\ess{\text{\rm{ess}}}
\def\beg{\begin} \def\beq{\begin{equation}}  \def\F{\scr F}
\def\Ric{\text{\rm{Ric}}} \def\Hess{\text{\rm{Hess}}}
\def\e{\text{\rm{e}}} \def\ua{\underline a} \def\OO{\Omega}  \def\oo{\omega}
 \def\tt{\tilde} \def\Ric{\text{\rm{Ric}}}
\def\cut{\text{\rm{cut}}} \def\P{\mathbb P} \def\ifn{I_n(f^{\bigotimes n})}
\def\C{\scr C}      \def\aaa{\mathbf{r}}     \def\r{r}
\def\gap{\text{\rm{gap}}} \def\prr{\pi_{{\bf m},\varrho}}  \def\r{\mathbf r}
\def\Z{\mathbb Z} \def\vrr{\varrho} \def\ll{\lambda}
\def\L{\scr L}\def\Tt{\tt} \def\TT{\tt}\def\II{\mathbb I}
\def\i{{\rm in}}\def\Sect{{\rm Sect}}  \def\H{\mathbb H}
\def\M{\scr M}\def\Q{\mathbb Q} \def\texto{\text{o}} \def\LL{\Lambda}
\def\Rank{{\rm Rank}} \def\B{\scr B} \def\i{{\rm i}} \def\HR{\hat{\R}^d}
\def\to{\rightarrow}\def\l{\ell}\def\iint{\int}
\def\EE{\scr E}\def\Cut{{\rm Cut}}
\def\A{\scr A} \def\Lip{{\rm Lip}}
\def\BB{\scr B}\def\Ent{{\rm Ent}}\def\L{\scr L}
\def\R{\mathbb R}  \def\ff{\frac} \def\ss{\sqrt} \def\B{\mathbf
B}
\def\N{\mathbb N} \def\kk{\kappa} \def\m{{\bf m}}
\def\dd{\delta} \def\DD{\Delta} \def\vv{\varepsilon} \def\rr{\rho}
\def\<{\langle} \def\>{\rangle} \def\GG{\Gamma} \def\gg{\gamma}
  \def\nn{\nabla} \def\pp{\partial} \def\E{\mathbb E}
\def\d{\text{\rm{d}}} \def\bb{\beta} \def\aa{\alpha} \def\D{\scr D}
  \def\si{\sigma} \def\ess{\text{\rm{ess}}}
\def\beg{\begin} \def\beq{\begin{equation}}  \def\F{\scr F}
\def\Ric{\text{\rm{Ric}}} \def\Hess{\text{\rm{Hess}}}
\def\e{\text{\rm{e}}} \def\ua{\underline a} \def\OO{\Omega}  \def\oo{\omega}
 \def\tt{\tilde} \def\Ric{\text{\rm{Ric}}}
\def\cut{\text{\rm{cut}}} \def\P{\mathbb P} \def\ifn{I_n(f^{\bigotimes n})}
\def\C{\scr C}      \def\aaa{\mathbf{r}}     \def\r{r}
\def\gap{\text{\rm{gap}}} \def\prr{\pi_{{\bf m},\varrho}}  \def\r{\mathbf r}
\def\Z{\mathbb Z} \def\vrr{\varrho} \def\ll{\lambda}
\def\L{\scr L}\def\Tt{\tt} \def\TT{\tt}\def\II{\mathbb I}
\def\i{{\rm in}}\def\Sect{{\rm Sect}}  \def\H{\mathbb H}
\def\M{\scr M}\def\Q{\mathbb Q} \def\texto{\text{o}} \def\LL{\Lambda}
\def\Rank{{\rm Rank}} \def\B{\scr B} \def\i{{\rm i}} \def\HR{\hat{\R}^d}
\def\to{\rightarrow}\def\l{\ell}\def\BB{\mathbb B}
\def\8{\infty}\def\I{1}\def\U{\scr U} \def\n{{\mathbf n}}\def\v{V}
\maketitle

\begin{abstract}  For a class of McKean-Vlasov stochastic differential equations with singular interactions, which include
the Coulomb/Riesz/Biot-Savart kernels as   typical examples (Examples 2.1 and 2.2), we derive  the well-posedness and regularity estimates by establishing the entropy-cost inequality. To measure   the singularity of interactions, we introduce a new probability distance
induced by local integrable functions, and estimate  this distance for the  time-marginal laws  of solutions by using  the Wasserstein distance of initial distributions.
A key point  of the study is to  characterize   the path space of  time-marginal distributions  for the solutions, by using  local hyperbound estimates on diffusion semigroups.
 \end{abstract} \noindent

\section{Introduction}
Let $\scr P$ be the set of all probability measures on $\R^d$ equipped with the weak topology. Consider
the following McKean-Vlasov stochastic differential equations on $\R^d$:
\beq\label{E0} \d X_t= b_t(X_t, \L_{X_t})\d t+  \si_t(X_t)\d W_t,\ \ t\ge 0,\end{equation}
where    $(W_t)_{t\ge 0}$ is an $m$-dimensional Brownian motion on a probability base (i.e.,  a complete filtered  probability space) $(\OO,\{\F_t\}_{t\ge 0},\F,\P)$, $\L_{X_t}$ is the distribution of $X_t$, and
$$ \si: [0,\infty)\times\R^d\to \R^d\otimes\R^m,\ \ \ b: [0,\infty)\times \R^d\times\tt{\scr P}\to \R^d$$
are measurable, where $\tt{\scr P}$ is a measurable subspace of $\scr P$ to be determined in terms of the singularity of $b_t(x,\cdot)$. When different probability spaces are concerned, we denote the distribution of $X_t$ under $\P$ by  $\L_{X_t|\P}$ to emphasize the underlying probability $\P$. To emphasize the distribution dependent property of \eqref{E0}, in the rest of this paper we call it distribution dependent stochastic differential equation (DDSDE).

Under  local integrability conditions on the time-spatial variables, as well as Lipschitz continuity of $b_t(x,\cdot)$ in Wasserstein or/and weighted variation distances,
the well-posedness, regularity estimates and ergodicity of \eqref{E0} have been extensively investigated,
see the recent monograph \cite{RW24} and references therein. There are also plentiful references
concerning other properties of this type SDEs, such as propagation of chaos and mean-field controls,
see for instance \cite{CD,JW} and references therein. 

In this paper, we aim to study the well-posedness and regularity estimates for \eqref{E0} with singular interactions, where the drift $b$ contains a term given by e.g.
\beq\label{0}b^{(0)}(x,\mu):= \int_{\R^d} {\bf K}(x,y)\mu(\d y),\ \ \ x\in\R^d,\ \mu\in\tt{\scr P}\end{equation}
for a   measurable map
$${\bf K}: \R^d\times\R^d\to\R^d$$ such that for each $x\in\R^d$, ${\bf K}(x,\cdot)$ is locally integrable with respect to the Lebesgue measure, and  $\tt {\scr P}$ is chosen such that the integral exists for $\mu\in \tt{\scr P}$.
Typical examples of ${\bf K}$ include the Coulomb/Newton, Riesz and Biot-Savart kernels, see \cite{L,S}:
  \beg{enumerate}
\item[(1)] {\bf Coulomb/Newton kernels.} Let $\oo_d$ be the volume of the unit ball in $\R^d$. The $d$-dimensional  Coulomb kernel
$${\bf K}_{C}(x,y):= \ff{x-y}{d\oo_d|x-y|^d},\ \ x\ne y$$
  describes electrostatic interactions between numerators; and the  Newton kernel ${\bf K}_{N}:=-{\bf K}_{C}$ reflects  gravitation interactions between bodies.
\item[(2)] {\bf Biot-Savart kernel.} Let $s_{d-1}$ be the area of $(d-1)$-dimensional unit sphere for $d\ge 2$, and let $z^\perp:=(-z_2,z_1)$ for $z=(z_1,z_2)\in\R^2$. The Biot-Savart kernel
$${\bf K}_{BS}(x,y):=\beg{cases} \ff{(x-y)^{\perp}}{2\pi|x-y|^2},\ &\text{if}\ d=2, \ x\ne y,\\
\ff{x-y}{s_{d-1}|x-y|^{d}},\ &\text{if}\ d\ge 3,\ x\ne y \end{cases}$$
describes interactions from incompressible fluids.
\item[(3)] {\bf Riesz  kernel.} For $0\ne \kk \in \R$ and $\bb\in (0, d)$, the Riesz  kernel
$${\bf K}_R(x,y):=\ff{\kk(x-y)}{|x-y|^{\bb+1}},\ \ x\ne y $$
   covers the Coulomb/Newton kernel and Boit-Savart kernel ($d\ge 3$), and has been applied in  solid state physics, ferrofluids and elasticity.
 \end{enumerate}

To characterize the singularity of $\mu\mapsto b^{(0)}(x,\mu)$ in $\eqref{0}$ with these singular kernels, we introduce below the new probability distance $\|\cdot\|_{k*}$ for $k\ge 1$ induced by $\tt L^k$-integrable functions.
Let    $\|\cdot\|_{L^k}$ be the $L^k$-norm with respect to the Lebesgue measure on $\R^d$, and denote
$$B(x,r):=\{y\in\R^d:\ |x-y|\le r\},\ \ \  (x,r)\in\R^d\times (0,\infty).$$
According to \cite{XXZZ}, $\tt L^k $ is the space of measurable functions $f$ on $\R^d$ such that
$$\|f\|_{\tt L^k}:=\sup_{x\in \R^d}\big\|1_{B(x,1)}f\big\|_{L^k}<\infty,\ \ \ k\in [1,\infty).$$
 Moreover, when $k=\infty$ we set
$$\|f\|_{\tt L^\infty}=\|f\|_{L^\infty}= \|f\|_\infty:=\sup_{x\in\R^d}|f(x)|.$$

 If  $|{\bf K}(x,y)|\le \ff c {|x-y|^\bb}$  for some constants $c>0$ and $\bb\in (0,d)$, which includes the above mentioned kernels as typical examples, then for   any
$k\in [1,\ff d\bb)$, we have
$$\sup_{x\in \R^d} \|{\bf K}(x,\cdot)\|_{\tt L^k}\le    \int_{B(0,1)}\ff c {|y|^{k\bb} }\d y =:K<\infty,$$
so that   the singular drift $b^{(0)}$ in \eqref{0} satisfies
$$|b^{(0)}(x,\mu)-b^{(0)}(x,\nu)|\le K \sup_{\|f\|_{\tt L^k}\le 1} |(\mu-\nu)(f)|,\ \ \mu,\nu\in \scr P_{k*},\ x\in\R^d,$$
where
\beq\label{PK*}   \scr P_{k*} :=\bigg\{\mu\in \scr P:\ \|\mu\|_{k*}:= \sup_{\|f\|_{\tt L^k}\le 1}\mu(|f|)<\infty\bigg\}.\end{equation}
Hence,  it is natural to study \eqref{E0} with such a singular interaction by using the $k*$-distance
\beq\label{PKD} \|\mu-\nu\|_{k*}:=\sup_{\|f\|_{\tt L^k}\le 1} |\mu(f)-\nu(f)|,\ \ \ \mu,\nu\in \scr P_{k*}.\end{equation}
Note that $k*$ here does not stand for the conjugate number $k^*:=\ff{k}{k-1}$, but refers to the dual norm
for measures induced by the $\tt L^k$ norm for functions.

For any  $k\in [1,\infty)$,    $(\scr P_{k*},\ \|\cdot\|_{k*})$ defined in \eqref{PK*} and \eqref{PKD} is a complete metric space, and the  Borel $\si$-field coincides with that induced by the weak topology, see Lemma \ref{L1} below.
When $k=\infty,$ we  set  $\scr P_{\infty*}:=\scr P$ and for any $\mu,\nu\in \scr P$,
\beg{align*}&\|\mu\|_{\infty*} :=\sup_{\|f\|_\infty\le 1} \mu(|f|) =1,\\
& \ \|\mu-\nu\|_{\infty*}= \|\mu-\nu\|_{var}:=\sup_{\|f\|_\infty\le 1} |\mu(f)-\nu(f)|.\end{align*}
So,  $(\scr P_{\infty*}, \|\cdot\|_{\infty*})= (\scr P,\|\cdot\|_{var})$  is complete as well.

It is clear that for constants $p\ge k\ge 1,$
$\|\cdot\|_{\tt L^k}\le \oo_d^{\ff{p-k}{pk}}\|\cdot\|_{\tt L^p}$, so that
$$\oo_d^{\ff{p-k}{pk}}\|\cdot\|_{k*}\ge \|\cdot\|_{p*},$$
hence the space $\scr P_{k*}$ is increasing in $k\ge 1$.

  To solve the SDE \eqref{E0} with the above mentioned singular interactions,  we consider solutions  satisfying $\L_{X_t}\in \tt{\scr P}:=\scr P_{k*}$ for some $k\in (1,\infty)$ such that $b^{(0)}(\cdot,\L_{X_t})$ is well-defined.
  To this end, for any $T\in (0,\infty)$ we shall introduce a path space $\C^T$ including weakly continuous maps from  $[0,T]$ to $\scr P$, such that for any $\mu=(\mu_t)_{t\in [0,T]}\in \C^T,$ the decoupled SDE
  \beq\label{DC} \d X_t^\mu= b_t(X_t^\mu,\mu_t)+\si_t(X_t^\mu)\d W_t,\ \ \ t\in [0,T],\ \L_{X_0^\mu}=\L_{X_0}\end{equation} with frozen distribution parameter $\mu$
  has a unique weak solution, and the map
$$\Phi: \mu\to \Phi\mu:=(\L_{X_t^\mu})_{t\in  [0,T]}$$
has a unique fixed point $\bar\mu$ in $\C^T$. If so, then $(X_t)_{t\in [0,T]}:=(X_t^{\bar\mu})_{t\in  [0,T]}$ is the unique weak solution of \eqref{E0} with $(\L_{X_t})_{t\in [0,T]}\in \C^T$.

Due to the regularization of noise, we may allow the initial distribution coming from a larger space  $\scr P_{p*}$ than $\scr P_{k*} $ for some $p>k$. In this case,  we should have   $\|\L_{X_t}\|_{k*}\to\infty $   as $t\to 0$
for $\L_{X_0}\in  \scr P_{p*}\setminus\scr P_{k*} .$ To describe this small time singularity, we recall the  local hyperbound estimate for a nice elliptic diffusion semigroup $P_t$ (see e.g. \cite{FYW23JDE}):
for any $T\in (0,\infty)$, there exists a constant $C(T)\in (0,\infty)$ such that
\beq\label{Hyp}\|P_t\|_{\tt L^k\to \tt L^p}:=\sup_{\|f\|_{\tt L^k}\le 1} \|P_tf\|_{\tt L^p} \le C(T) t^{-\ff{d(p-k)}{2pk}},\ \ t\in (0,T],\ \infty\ge p\ge k\geq 1, \end{equation}
where $\ff{d(p-k)}{2pk}:= \ff d{2k}$ when $p=\infty$.   If this estimate holds for the diffusion semigroup associated with \eqref{DC}, then    for any initial distribution $\gg:=\L_{X_0}\in \scr P_{p*}$,   the time-marginal distribution
$(\L_{X_t})_{t\in  [0,T]}$ of solution  to \eqref{DC}  up to time $T$     belongs  to the  path space
 \beq\label{CPK}\C_{p,k}^{T}:=\Big\{\mu\in C^w([0,T];\scr P):\ \rr_{T}^{p,k}(\mu):=\sup_{t\in (0,T]} t^{\ff{d(p-k)}{2pk}}\|\mu_t\|_{k*}<\infty\Big\},\end{equation}
where $C^w([0,T];\scr P)$ is the set of all weakly continuous maps from $[0,T]$ to $\scr P$. This leads to the following  notion of the maximal $\C_{p,k}$-solution for \eqref{E0}, where the life time is the smallest time $\tau\in (0,\infty)$ such that
$\limsup_{t\uparrow\tau}\|\L_{X_t}\|_{k*}=\infty$, and we denote $\tau=\infty$ if such a finite time does not exist. Since  $\L_{X_t}$ is  deterministic, so is the  life time $\tau$.

\beg{defn}[Maximal strong $\C_{p,k}$-solution]  Let $k\in [1,\infty]$ and  $p\in [k,\infty]$.
  We call $(X_t)_{t\in [0,\tau)}$ a  maximal strong  $\C_{p,k}$-solution of \eqref{E0} with life time $\tau$, if it is an adapted continuous process on $\R^d$ such that the following conditions hold.
\beg{enumerate}
\item[(1)] The initial distribution $\L_{X_0}\in \scr P_{p*}$, $\tau\in (0,\infty],$      and
$$\limsup_{t\uparrow \tau} \|\L_{X_t}\|_{k*}=\infty \ \ \text{if}\ \tau<\infty.$$
\item[(2)] For any $T\in  (0,\tau)$, $(X_t)_{t\in [0,T]}$ is a strong $\C_{p,k}$-solution of \eqref{E0} up to time $T$, i.e.
$$(\L_{X_t})_{t\in [0,T]}  \in \C_{p,k}^T,\ \
 \E \int_0^T \big[| b_s(X_s,\L_{X_s})|+ \|\si_s(X_s)\|^2\big]\d s<\infty,$$   and $\P$-a.s.
$$X_t= X_0+\int_0^t b_s(X_s,\L_{X_s})\d s+ \int_0^t \si_s(X_s)\d W_s,\ \ t\in [0,T].$$
\end{enumerate}
When  $\tau=\infty$, we call $(X_t)_{t\ge 0}$  a  global strong $\C_{p,k}$-solution of \eqref{E0}. \end{defn}

\beg{defn}[Maximal weak  $\C_{p,k}$-solution]  Let $k\in [1,\infty]$, $p\in [k,\infty]$ and $\gg\in \scr P_{p*}$.
\beg{enumerate}
\item[(1)]    A couple $(X_t,W_t)_{t\in [0,\tau)}$ is called a maximal weak $\C_{p,k}$-solution of \eqref{E0} with initial distribution $\gg$,
 if there exists a probability base $(\OO,\{\F_t\}_{t\in [0,\tau)},\F,\P)$
such that $(W_t)_{t\in [0,\tau)}$ is an $m$-dimensional Brownian motion, $\L_{X_0}=\gg$ and $(X_t)_{t\in [0,\tau)}$ is a maximal strong $\C_{p,k}$-solution of
\eqref{E0}. In this case, for any $ T\in (0,\tau)$, $(X_t,W_t)_{t\in [0,T]}$ is called a   weak $\C_{p,k}$-solution of \eqref{E0} up to time $T$.
\item[(2)]  If \eqref{E0} has a maximal weak $\C_{p,k}$-solution with initial distribution $\gg$,   and   any two maximal weak $\C_{p,k}$-solutions with initial distribution $\gg$ have common life time and distribution,
   then we say that \eqref{E0} with initial distribution $\gg$ has a unique maximal weak $\C_{p,k}$-solution. In this case, we denote the life time by $\tau(\gg)$ and set
$$P_t^*\gg:=\L_{X_t},\ \ t\in [0,\tau(\gg)).$$
 \end{enumerate}\end{defn}

Note that for any $T\in (0,\infty)$ and $\mu\in \C_{p,k}^{T}$,
$$\|\mu_t\|_{k*}\le c t^{-\ff{d(p-k)}{2pk}},\ \ \ t\in (0,T]$$
holds for some constant $c\in (0,\infty)$. So, to ensure  $\int _0^{T} \|\mu_t\|_{k*}^{2}\d t<\infty$, which is essential to apply Girsanov's theorem with
drift having   linear growth  in  $\|\mu_t\|_{k*}$, we   need $\ff{ d(p-k)}{pk}<1$, i.e. $(p,k)$ belongs to the class
 \beq\label{DK} \D :=\Big\{(p,k):\ 1\le k\le p\le \infty,  \ \ff 1 k-\ff 1 d<\ff 1 p\Big\}.\end{equation}

To  cover more general models, besides a drift term $b^{(0)}$ as in  \eqref{0} with singular interaction, we also consider two additional drift terms: the regular term $b_t^{(1)}$ is Lipschitz continuous on $\R^d\times \scr P_{k*}$,  and  the singular term $\sum_{i=2}^{l_0} b^{(i)} $ for some $2\le l_0\in \mathbb N$  satisfying  time-spatial local integrability conditions. So, the drift $b_t$ is decomposed as
\beq\label{B} b_t(x,\mu)= b_t^{(0)}(x,\mu)+ b_t^{(1)}(x,\mu)+ \sum_{i=2}^{l_0} b_t^{(i)}(x,\mu).\end{equation}

 In Section 2, we state the main results of the paper concerning the well-posedness  (i.e. existence and uniqueness) and regularity estimates  for   the maximal strong/weak $\C_{p,k}$-solutions  of \eqref{E0},
 which are illustrated by    typical  examples  of  the above mentioned singular kernels. The proofs of these results will be addressed in Section  3 and  Section 4,
 with helps of  preliminary results introduced in Section 5, where some existing results on singular SDEs are extended to the case with several singular drifts.

\section{Main results and examples}

As explained above, we shall use some $k*$-distance to measure the singularity of interactions. To characterize the time-spatial singularity, we  recall the family of locally integrable functions
introduced in \cite{XXZZ}.

For any $p,q\in [1,\infty)$ and $0\le s<t<\infty$, let $\tt L_{q}^{p}(s,t)$
be the set of measurable functions $f: [s,t]\times\R^d\to \R$ such that
$$\|f\|_{\tt L_{q}^{p}(s,t)}:=\sup_{x\in\R^d} \bigg(\int_s^t\|1_{B(x,1)}f_r\|_{L^{p}}^{q}\d r\bigg)^{\ff 1 {q}}<\infty.$$    Simply denote
$\tt L_{q}^{p}(t):= \tt L_q^p(0,t), \|\cdot\|_{\tt L_{q}^{p}(t)}:= \|\cdot\|_{\tt L_{q}^{p}(0,t)}. $

We will take  $(p,q)$  from the following class
$$\scr K:= \Big\{(p,q)\in (2,\infty):\ \ff d {p}+\ff 2 {q}<1\Big\}.$$
 We make the following assumptions.

\beg{enumerate} \item[{\bf (A)}] Let $(p,k)\in\D$ defined in \eqref{DK},     $b^{(i)} (0\le i\le l_0)$ be in \eqref{B}.  For any $T>0, (t,x)\in[0,T]\times\R^d$ and $\mu\in \C_{p,k}^T, $ denote
$$a_t(x):=(\si_t\si_t^*)(x),\ \ b^{i,\mu}_t(x):= b^{(i)}_t(x,\mu_t),\ \     \  2\le i\le l_0.$$
\item[$(A_1)$] For any $T\in (0,\infty),$ there exist  $K\in (0,\infty)$,   $\aa\in (0,1]$  and
$\{(p_i',q_i'): 2\le i\le l_0\}\subset \scr K$   such that
for any  $t\in [0,T], x,y\in\R^d,$   $\nu,\tt\nu\in \scr P_{k*} $ and $\mu\in \C_{p,k}^T$,
$$ |b_t^{(0)}(x,\nu) | \le K \|\nu\|_{k*},\ \ \ \ |b_t(x,\nu)-b_t(x,\tt \nu)|\le K\|\nu-\tt\nu\|_{k*},$$
$$ b_t^{(1)}(0,\mu)=0,\ \ \ \   \
    |b_t^{(1)}(x,\nu)-b_t^{(1)}(y,\tt\nu)|\le K(|x-y|+\|\nu-\tt\nu\|_{k*}),$$
$$\|a\|_\infty+\|a^{-1}\|_\infty+\sup_{2\le i\le l_0}\|b^{i,\mu} \|_{\tt L_{q_i'}^{p_i'}(T)}\le K,\ \ \     |a_t(x)-a_t(y)|\le K|x-y|^\aa.$$
\item[$(A_2)$] For any $T\in (0,\infty)$, $a_t(x)$ is weakly differentiable in $x\in\R^d$ for a.e. $t\in [0,T]$,  and  there exist  finite many  $(p_i,q_i)\in \scr K$ and   $1\le f_i\in \tt L_{q_i}^{p_i}(T) $  for   $1\le i\le \ell$,   such that $$\|\nn a\|\le \sum_{i=1}^\ell f_i.$$
\end{enumerate}

\beg{thm}\label{T0}  Assume $(A_1)$ and let $b$ be in $\eqref{B}$.  Then the following assertions hold.
 \beg{enumerate} \item[$(1)$]   For any initial distribution $\gg\in \scr P_{p*}$,   $\eqref{E0}$   has a unique maximal weak $\C_{p,k}$-solution with life time
 $\tau(\gg)\in (0,\infty].$
  \item[$(2)$]  For any $n\in \mathbb N,$ there exist constants $\bb_0(n)\in (0,1]$ and $\bb_1(n)\in [1,\infty)$ such that
\beq\label{TT0} \tau(\gg)>\tau_n(\gg):=\beg{cases}n,\ &\text{if} \   p=\infty \ \text{or }\  b^{(0)}=0,\\
 \bb_0(n)\|\gg\|_{p*}^{-1/\theta},\ &\text{otherwise,}\end{cases}\end{equation}  where  $\theta:=\ff 1 2 -\ff{d(p-k)}{2pk}>0,$ and
\beq\label{EST}  \sup_{t\in (0,\tau_n(\gg)]} t^{\ff{d(p-k)}{2pk}}\|P_t^*\gg\|_{k*}\le \bb_1(n)\|\gg\|_{p*},\ \ \gg\in \scr P_{p*}.\end{equation}
 \item[$(3)$]  If $\tau(\gg)<\infty$, then
\beq\label{L} \liminf_{t\uparrow\tau(\gg)} \big(\tau(\gg)-t\big)^{\theta} \|P_t^*\gg\|_{p*}>0,\end{equation}
\beq\label{L'} \int_r^{\tau(\gg)}  \|P_t^*\gg\|_{k*}^2   \d t=\infty,\ \ \ r\in [0,\tau(\gg)).  \end{equation}
 \item[$(4)$] If   $(A_2)$ holds,   then for any $\F_0$-measurable initial value $X_0$ with $\gg:=\L_{X_0}\in \scr P_{p*}$, the SDE $\eqref{E0}$ has a unique maximal strong $\C_{p,k}$-solution.
   Moreover,  there exists an  increasing function   $C_\gg: [1,\infty)\times   (0,\tau(\gg)) \to (0,\infty)$ such that
\beq\label{NES} \E\bigg[\sup_{s\in [0,t]} |X_s|^n\bigg|\F_0\bigg]\le C_\gg(n,t) (1+|X_0|^n),\ \ n\in [1,\infty),\ t\in (0,\tau(\gg)).\end{equation}
If either $p=\infty$ or $b^{(0)}=0$, then $\tau(\gg)=\infty$ and   $C_\gg(n,t)=C(n,t)$ is independent of $\gg\in \scr P_{p*}$.
 \end{enumerate}
\end{thm}

\paragraph{Remark 2.1.}  Theorem \ref{T0}(3) shows that   the blowup in the larger   $k*$-distance  is equivalent to that in the smaller  $p*$-distance for the maximal $\C_{p,k}$-solution, where
      \eqref{L}  is in the same spirit of    Leray's   blowup criterion \cite{L} for 3D Navier-Stokes equation, and \eqref{L'} implies that for any constant $\kk>\ff 1 2$,
      $$\limsup_{t\uparrow \tau(\gg)} \|P_t^*\gg\|_{k*}\ss{\tau(\gg)-t} \Big(\log \Big[1+  \big(\tau(\gg)-t\big)^{-1}\Big]\Big)^\kk =\infty\ \text{if}\ \tau(\gg)<\infty.$$
We would like to compare Theorem \ref{T0} with some existing results  for SDEs with singular interactions.
\beg{enumerate}
\item[(1)]   Let $\dd_x$ denote the Dirac measure at $x\in \R^d$.  When  $a:=\si\si^*$ satisfies  $(A_1)$,  and $b$  satisfies
\beq\label{*}\|b_\cdot(\cdot,\dd_0)\|_{\tt L_{q_0}^{p_0}(T)}\le K,\ \   \|b_\cdot(\cdot, \gg)-b_\cdot(\cdot, \tt\gg)\|_{\tt L_{q_0}^{p_0}(T)}\le K\|\gg-\tt\gg\|_{var}\end{equation}
for some constants $T,K\in (0,\infty)$ and $(p_0,q_0)\in \scr K$,    the weak well-posedness of \eqref{E0} up to time $T$  has been presented in  \cite[Theorem 1.1]{21RZ} and \cite[Proposition 1.2]{Zhao}.  It is in particular the case when
\beq\label{KL} |{\bf K}(x,y)|\sim  \ff 1 {|x-y|^\bb}\,\ \ x\ne y \end{equation}   for\ some $\bb\in (0,1).$
Since   $(A_1)$   uses larger probability distance $\|\cdot\|_{k*}$ instead of $\|\cdot\|_{var}= \|\cdot\|_{\infty*},$ Theorem \ref{T0}  applies to  examples which do not satisfy   \eqref{*}. For instance,
when $b=b^{(0)}$ defined in \eqref{0} for the kernel in \eqref{KL} with $\bb\in [1,d)$ for $d\ge 2$, and \eqref{*} does not hold but $(A_1)$ does when
 $$ {\bf K}(x,y)= \ff 1 {|x-y|^\bb}+\ff 1 {|y|^\bb},\ \ x\ne y $$ for\ some $\bb\in (0,d).$
\item[(2)]
When $a$ is the identity matrix $I_{d\times d}$,  the SDE \eqref{E0} with drift $b=b^{(0)}$ given by \eqref{0} has been investigated in many papers, in particular for  ${\bf K}={\bf K}_{BS},$  see
\cite{BRZ, B94,GG05,Kato94} and references within. For ${\bf K}$ in \eqref{KL} with  some constants $c\in (0,\infty)$ and  $\bb\in [1,d)$, the weak well-posedness of \eqref{E0} up to a deterministic time  $T\sim \|\ell_\gamma\|_{\infty}^{-2}$ has been derived in  \cite[Theorem 2]{LQ}, see also \cite[Theorem 1.1]{CQ} and \cite[Proposition 3.1]{RS} for the locally weak well-posedness of the associated non-linear Fokker-Planck equation, where $\ell_\gg:=\ff{\d\gg}{\d x}$ is not necessarily bounded.
Note that in this case   {\bf (A)} holds for  any  $k\in (1,\ff 3{2})$  and $p\in   [k, \ff{3k}{3-k})$,  so that Theorem \ref{T0} ensures the weak and strong well-posedness for $\C_{p,k}$-solutions of  \eqref{E0}  for any initial distribution with $\|\gamma\|_{p*}<\infty$ up to a time $ T \sim \|\gamma\|_{p*}^{-1/\theta}$.
\item[(3)]    We will show in Corollary \ref{C01} that \eqref{E0} is globally well-posed for $\C_{p,k}$-solution when the associated Fokker-Planck equation is well-posed for solutions with bounded densities, which is, in particular,  the case when ${\bf K}$  is the 2D  Biot-Savart kernel.
   \end{enumerate}

As a consequence of  Theorem \ref{T0},    we have  the following  criteria  on the global well-posedness  of  \eqref{E0} by using the associated nonlinear  Fokker-Planck equation:
\beq\label{NFK} \pp_t \mu_t= L_{\mu_t}^*\mu_t,\ \ L_{\mu_t}:= \ff 1 2 {\rm tr}(a_t\nn^2)+ b_t(\cdot,\mu_t)\cdot\nn,\ \ t\ge s.\end{equation}
 A solution of this PDE   is a weak continuous map $\mu_\cdot: [s,\infty)\to \scr P $ such that
$$\mu_t(f_t)=\mu_s(f_s) +\int_s^t L_{\mu_r} f_r \d\mu_r,\ \ \ f\in C_0^\infty([s,\infty)\times \R^d),\ t\in [s,\infty).$$

\beg{cor}  \label{C01} Assume  $(A_1)$.  Let $b$ be in $\eqref{B}$ with $b^{(i)}=0$ for $2\le i\le l_0$,   and let  $\gg\in \scr P_{p*}$ such that
$\gg(|\cdot|)<\infty$ when $b^{(1)}\ne 0$.
   If    there exists $s\in (0,\tau(\gg))$ such that for any   $\mu_s\in \scr P $ with $\|\ell_{\mu_s}\|_\infty<\infty$,   the PDE $\eqref{NFK}$ for $t\ge s$ has a global solution $(\mu_t)_{t\ge s}$ with
\beq\label{*E}  \sup_{t\in [s,T]  } \|b_t(\cdot,\mu_t)\|_\infty <\infty,\ \ T\in [s,\infty),\end{equation}
then \eqref{E0} has a unique  global weak $\C_{p,k}$-solution (i.e. $\tau(\gg)=\infty$), and
 \beq\label{ER}  \sup_{t\in (0,T]}  t^{\ff{d(p-q)}{2qp}}  \|P_t^\ast\gamma\|_{q\ast}<\infty,\ \ q\in [1,p],\ T\in (0,\infty).\end{equation}
 If moreover $(A_2)$ holds, then for any initial value $X_0$ with $\L_{X_0}=\gg$,  \eqref{E0}   has a unique global strong $\C_{p,k} $-solution.
 \end{cor}

 By combining Corollary \ref{C01} with the well-posedness of   2D Navier-Stokes which has been well-studied  in the literature of PDEs,
 we present below an example ensuring  the global well-posedness of  strongly $\C_{p,k}$-solution  for  the DDSDE \eqref{E0} with interaction given by the 2D Biot-Savart kernel.
 This will enable us to establish the entropy-cost inequality in Example 2.3(3) below, which is new  from both  literatures of PDEs and SDEs.

 \paragraph{Example 2.1.}  Let $d=2, \si=\kk I_{2\times 2}$ for some constant $\kk\in (0,\infty)$, and
 $$b_t(x,\mu):= \int_{\R^d} {\bf K}(x-y) \mu(\d y)$$ for the Biot-Savart kernel
$ {\bf K}(x):= \ff {(-x_2,x_1)} {2\pi |x|^2},\   x=(x_1,x_2)\in\R^2.$
 Then for any $k\in (1,2)$,   $p\in  [k,\ff{2k}{2-k})$  and $\gg\in \scr P_{p*}$, the SDE \eqref{E0} has a unique global strong $\C_{p,k}$-solution, and \eqref{ER} holds.

 \beg{proof}
 For a fixed $s\in (0,1\land \tau(\gg))$,  consider   the
  2D   vorticity   equation
 \beq\label{VNS} \pp_t v_t= \ff {\kk^2} 2  \DD v_t- (u_t\cdot\nn)v_t,\ \ u_t(x):= \int_{\R^d} {\bf K}(x-y)v_t(y)\d y,\ \ t\in [s,\infty)   \end{equation}
This equation is equivalent to  \eqref{NFK} for $b_t=u_t$.
 By \cite[Theorem 4.3]{GMO},  for any  probability density $\|v_s\|_\infty<\infty$, \eqref{VNS} has a  unique global  solution  with
 $$\sup_{t\in [s,T]}  \|v_t\|_\infty<\infty,\ \ T\in (s,\infty).$$   Then    $b=b^{(0)}:=u$ and $\mu_t(\d x):= v_t(x)\d x$ satisfy
\beg{align*}&\sup_{t\in [s,T]  } \|b_t(\cdot,\mu_t)\|_\infty \le  1+ \sup_{t\in [s,T],x\in\R^2  } \int_{B(x,1)} \ff {v_t(y)} {|y-x|}  \d y \\
&\le 1 +  \Big(  \sup_{t\in [s,T]}  \|v_t\|_\infty\Big)\int_{B(0,1)} \ff{\d y}{|y|}  <\infty. \end{align*}   So,
 \eqref{*E} holds and  the desired assertion follows from Corollary \ref{C01}.
 \end{proof}

Having  the maximal weak well-posedness for the $\C_{p,k}$-solution of \eqref{E0}, our main concern   is to    study the regularity of  the map
$$\scr P_{p*}\ni \gg\mapsto P_t^*\gg\in \scr P_{k*}$$ for $t\in (0,\tau(\gg))$
by estimating the $k*$-distance
  $\|P_t^*\gg-P_t^*\tt\gg\|_{k*}$ and the relative entropy $\Ent(P_t^*\gg|P_t^*\tt\gg),$
  using the Wasserstein distances $\W_q(\gg,\tt\gg)$ for some
$q\ge 1.$    Recall that for any $\gg,\tt \gg\in \scr P,$
$$\Ent(\gg|\tt\gg):= \beg{cases} \gg\big(\log\ff{\d\gg}{\d\tt\gg}\big),\ &\text{if}\ \ff{\d\gg}{\d\tt\gg}\ \text{exists},\\
\infty,\ &\text{otherwise},\end{cases}$$ and for any constant $q\in [1,\infty)$,
$$\W_q(\gg,\tt\gg):=\inf_{\pi\in \C(\gg,\tt\gg)}\bigg(\int_{\R^d\times\R^d} |x-y|^q\pi(\d x,\d y)\bigg)^{\ff 1 q},$$   where $\C(\gg,\tt\gg)$ is the set of all couplings for $\gg$ and $\tt\gg$.
The estimates will depend on
\beq\label{1ga} \kk_t(\gg):= 1_{\{\|b^{(0)}\|_\infty>0\}}\Big(\|\gg\|_{p*}\lor \sup_{s\in (0,t]} s^{\ff{d(p-k)}{2pk}} \|P_s^*\gg\|_{k*}\Big),\ \ \ t\in (0,\tau(\gg)),\ \gg\in \scr P_{p*}.\end{equation}
By \eqref{EST}, $\kk_t(\gg)\le \bb_1(n)\|\gg\|_{p*}$ for $t\le \tau_n(\gg)$.

Recall that $\theta:=\ff 1 2-\ff{d(p-k)}{2pk}>0.$ For any $\gg,\tt\gg\in \scr P_{p*}$, $ t\in (0,\tau(\gg)\land\tau(\tt\gg))$ and  increasing function $\bb: (0,\infty)\to (0,\infty)$,  let
 \beq\label{Ka}  K_{t,\bb}^{p,k}(\gg,\tt\gg):=
 \exp\Big[ \bb_t  \e^{\bb_t(t\kk_t(\gg)^{1/\theta}+ t\kk_t(\tt\gg)^{1/\theta})}\Big].\end{equation}
  Moreover,   for any $\theta'\in (0,\theta),$ let
\beq\label{ST}s_t(\theta',\gg):=   \beg{cases} t\land [\kk_t(\gg)^{-1/\theta'}],\ &\text{if} \ \|b^{(0)}\|_\infty>0,\\
 t,\ &\text{if}\ b^{(0)}\equiv 0.\end{cases}\end{equation}

 \beg{thm}\label{T02}  Let $b$ be in $\eqref{B}$ such that   {\bf (A)}  holds. Assume $k\in(\frac{d}{2},\infty]$.
 Then   for any  $q\in [1,\infty)$ such that  $(\frac{pq}{q-1},k)\in \D$,   where $\frac{pq}{q-1}:=\infty$ if $q=1,$ the following assertions hold for some increasing $\bb: [0,\infty)\to (0,\infty),$
  all $\gg,\tt\gg\in \scr P_{p*}$ and any $t\in (0,\tau(\gg)\land \tau(\tt\gg)).$
 \beg{enumerate} \item[$(1)$]  We have \beq\label{ES5} \beg{split}   \|P_t^*\gg-P_t^*\tt\gg\|_{k*}
 &\le (\|\gg\|_{p*}+\|\tt\gg\|_{p*})^{\ff{q-1}q} K_{t,\bb}^{p,k}(\gg,\tt\gg) t^{-\left((\ff 1 2+\ff{d(pq-(q-1)k)}{2pqk})\vee\frac{d}{2k}\right)}    \W_q(\gg,\tt\gg).  \end{split}\end{equation}
If either $p=\infty$ or $b^{(0)}=0,$  then
\beq\label{ES5'}   \|P_t^*\gg-P_t^*\tt\gg\|_{k*} \le \bb_t (\|\gg\|_{p*}+\|\tt\gg\|_{p*})^{\ff{q-1}q}   t^{-\left((\ff 1 2+\ff{d(pq-(q-1)k)}{2pqk})\vee\frac{d}{2k}\right)}\W_q(\gg,\tt\gg). \end{equation}
\item[$(2)$] For any $\theta'\in (0,\theta),$
\beq\label{ES6} \beg{split}  \Ent(P_t^*\gg|P_t^*\tt\gg)
  &\le \bb_t(\|\gg\|_{p*}+\|\tt\gg\|_{p*})^{\ff{2(q-1)}q}\\
&\quad\times   \bigg(\ff{\W_2(\gg,\tt\gg)^2}{s_t(\theta',\gg)}
+\ff  {K_{t,\bb}^{p,k}(\gg,\tt\gg)^2\W_q(\gg,\tt\gg)^2}{[s_t(\theta',\gg) \land s_t(\theta',\tt\gg)]^{\left(\left(1+\ff{d(qp-(q-1)k)}{pqk}\right)\vee\frac{d}{k}\right)-1}}\bigg). \end{split}\end{equation}
 In particular, if   $p=\infty$, then
 \beq\label{ES7}  \Ent(P_t^*\gg|P_t^*\tt\gg)
\le   \ff {\bb_t}t  \W_2(\gg,\tt\gg)^2,\ \  t> 0,  \end{equation}
while for $b^{(0)}=0$ and  $p<\infty$,
\beq\label{ES8}  \beg{split}  \Ent(P_t^*\gg|P_t^*\tt\gg)
  \le \  & \bb_t (\|\gg\|_{p*}+\|\tt\gg\|_{p*})^{\ff{2(q-1)}{q}}   \\
  &\times \bigg(\ff{ \W_2(\gg,\tt\gg)^2}{  t}  + \ff{ \W_q(\gg,\tt\gg)^2}{  t^{\left(\left(1+\ff{d(qp-(q-1)k)}{pqk}\right)\vee\frac{d}{k}\right)-1}}  \bigg),\ \ t> 0.\end{split}\end{equation}
\end{enumerate}
\end{thm}

\paragraph{Remark 2.2.} Since  $\|\cdot\|_{k*}$ is essentially larger than $\|\cdot\|_{var}$, we see that \eqref{ES5} is stronger than  the same type  estimates on $\|P_t^*\gg- P_t^*\tt\gg\|_{var}.$
The estimate   \eqref{ES6} is called the entropy-cost inequality or the log-Harnack inequality, which has been established for various models including SDEs, SPDEs and McKean-Vlasov SDEs, see for instance \cite{R23, RW21, Wbook, RW24} and references therein. This type estimate has been derived in \cite{HW*}   for $\ff 1 2$-Dini  interactions, see also   \cite{HRWJDE} for  the case with   distribution dependent noise, where
$$|b_t(x,\mu)-b_t(x,\nu)|\le K(\W_q+\W_\psi)(\mu,\nu)$$
holds for some constant $K\in (0,\infty)$ and the Wasserstein distance
$$\W_\psi(\mu,\nu):=\sup\bigg\{|\mu(f)-\nu(f)|:\ \sup_{x\ne y}\ff{|f(x)-f(y)|}{\psi(|x-y|)}\le 1\bigg\}$$
induced by an increasing concave function $\psi$ with $\psi(0)=0$ and $\int_0^t \ff{\psi(s)^2}{s}\d s<\infty$, i.e. $\psi^2$ is a Dini function so that $\W_\psi$ describes  $\ff 1 2$-Dini  interaction kernels.
However, when the interaction is singular of type   \eqref{0} with only  locally integrable kernels, the log-Harnack inequality is unknown until the present work.

\

To illustrate Theorem \ref{T02}, we present below an   example where the interaction is general enough to cover  the Coulomb/Riesz/Biot-Savart  kernels.

\paragraph{Example 2.2.} Let $b^{(1)}$, $b^{(i)} (2\le i\le l_0)$ and $a:=\si\si^*$ satisfy the corresponding conditions in {\bf (A)}, and let $b^{(0)}$ be in  \eqref{0}
   such that
$$|{\bf K}(x,y)|\le \ff{c }{|x-y|^\beta}+\sum_{i=1}^l \ff c{|y-x_i|^\bb},\ \ y\notin\{x, x_i:1\le i\le l\} $$ holds for some constants $c \in (0,\infty), \beta\in (0,d), l\in\mathbb N$ and $\{x_i:1\le i\le l\}\subset \R^d$.
Then all assertions in Theorem  \ref{T0} and Theorem \ref{T02}    hold for any $k\in \big(1,\ff d {\beta}\big), $ $p\in [k,\infty]$ and $q\in [1,\infty)$ such that $(\frac{pq}{q-1},k)\in \D$, i.e.
$\ff 1 k- \ff 1 d <\ff {q-1} {pq}.$
In particular:
\beg{enumerate}\item[$(1)$] If $\beta<1,$ then we may take $k\in \big(d,\ff d\beta\big)$ and $p=\infty$ such that \eqref{ES7} holds.
\item[(2)] When ${\bf K}$ is one of the Coulomb/Biot-Savart  kernels for $d\ge 2$,  all  assertions in Theorem   \ref{T02}, except \eqref{ES7} and \eqref{ES8},    hold for
$$k\in \Big(1,\ff{d}{d-1}\Big),\ \ \  p\in \Big[k, \ff{dk}{d-k}\Big),
 \ \ \ q\in \Big(\ff{dk}{dk-p(d-k)},\infty\Big).$$
\item[(3)]  In Example 2.1 where  {\bf K} is the 2D Biot-Savart kernel, Theorem \ref{T02} applies to  $k\in (1,2), p\in \big[k,\ff{2k}{2-k}\big)$ and $q\in \big(\ff{2k}{2k-p(2-k)},\infty\big)$, for
$\tau(\gg)=\infty$.  \end{enumerate}

\section{Proofs of Theorem \ref{T0} and Corollary \ref{C01}}

Given initial distribution $\gg\in \scr P_{p*}$ and  $T\in (0,\infty)$,    let $\C_{p,k}^{T}$ be in
\eqref{CPK} and
$$\C_{p,k}^{\gg,T}:=\big\{\mu\in \C_{p,k}^{T}:\ \mu_0=\gg\big\}.$$
The existence and uniqueness of (weak) $\C_{p,k}$-solution  of \eqref{E0} with $\L_{X_0}=\gg\in \scr P_{p*}$
up to time $T$ holds, where $T$ may depend on $\gg$, if we could verify the following assertions:
\beg{enumerate} \item[(i)] The metric space $(\C_{p,k}^{\gg,T}, \rr_{T}^{p,k})$ is complete for $\rr_{T}^{p,k}$ defined by
$$\rr_{T}^{p,k}(\mu,\nu):=\sup_{t\in (0,T]} t^{\ff{d(p-k)}{2pk}} \|\mu_t-\nu_t\|_{k*},\ \ \mu,\nu\in \C_{p,k}^{\gg,T}.$$
\item[(ii)] For any $\mu\in \C_{p,k}^{\gg,T}$, the SDE
\beq\label{EM} \d X_t^\mu=b_t(X_t^\mu,\mu_t)\d t+\si_t(X_t^\mu)\d W_t,\ \ t\in [0,T]\end{equation}
has a unique weak solution with initial distribution $\gg$ such that the element
\beq\label{PHG} \Phi^\gg\mu= (\Phi^\gg_t \mu)_{t\in [0,T]} :=(\L_{X_t^\mu})_{t\in [0,T]}\in \C_{p,k}^{\gg,T}.\end{equation}
\item[(iii)]  The map $\Phi^\gg: \C_{p,k}^{\gg,T}\to \C_{p,k}^{\gg,T}$ has a unique fixed point.
\end{enumerate}

Once these three items  are confirmed, letting $\mu$ be the unique fixed point of $\Phi^\gg$ in $\C_{p,k}^{\gg,T}$,
we see that $(X_t^\mu,W_t)_{t\in [0,T]}$ becomes the unique weak $\C_{p,k}$-solution of \eqref{E0} up to time $T$, and if \eqref{EM} has a unique strong solution with initial value $X_0$ such that
$\L_{X_0}=\gg$, then $(X_t^\mu)_{t\in [0,T]}$ is also the unique strong $\C_{p,k}$-solution of \eqref{E0} up to time $T$.
To verify the above assertions, we present below some lemmas.

\beg{lem}\label{L1} Let $k\in [1,\infty]$, $p\in   [k,\infty]$, $\ll\in [0,\infty)$ and $T\in (0,\infty)$.
Then the following assertions hold.
\beg{enumerate} \item[$(1)$]  The space $(\scr P_{k*},\ \|\cdot\|_{k*})$ defined in \eqref{PK*} and \eqref{PKD} is   complete, and the Borel $\si$-field coincides with that induced by the weak topology.
\item[$(2)$]
 The space $(\C_{p,k}^{\gg,T},\rr_{T,\ll}^{p,k})$ is complete, where
   $$\rr_{T,\ll}^{p,k}(\mu,\nu):=\sup_{t\in (0,T]} \e^{-\ll t}t^{\ff{d(p-k)}{2pk}} \|\mu_t-\nu_t\|_{k*}.$$
\end{enumerate}  \end{lem}

\beg{proof} (1)
For any $r\ge  \ss d,$ we find a constant $c(r)\in \mathbb N$ such that
each $B(x,1)$ is covered by $c(r)$ many sets in $\{B(z,r):z\in \Z^d\}$, while every $B(z,r)$ is covered by
$c(r)$ many sets in $\{B(x,1):x\in\R^d\}$. Hence,
\begin{align}\label{dfk}c(r)^{-1} \sup_{z\in\Z^d}\|f1_{B(z,r)}\|_{L^k} \le
  \|f\|_{\tt L^k} \le  c(r) \sup_{z\in \Z^d} \|1_{B(z,r)} f\|_{L^k},
\end{align}
 So,
  $\mu\in \scr P_{k*}$  implies that $\ell_\mu:=\ff{\d\mu}{\d x}$ exits, and
 \beq\label{SJ0} c(r)^{-1} \sum_{z\in\Z^d} \|\ell_\mu 1_{B(z,r)}\|_{L^{\ff k{k-1}}}\le \|\mu\|_{k*}\le c(r)\sum_{z\in\Z^d} \|\ell_\mu 1_{B(z,r)}\|_{L^{\ff k{k-1}}}.\end{equation}
 Indeed, by $\cup_{z\in\Z^d}B(z,r)=\R^d$
 and  noting that \eqref{dfk} implies
 $$\sup_{\|f\|_{\tt L^k}\le 1}\|f1_{B(z,r)}\|_{L^{k}}\le c(r),\ \ z\in\Z^d,$$
 we derive
 $$\|\mu\|_{k*}:= \sup_{\|f\|_{\tt L^k}\le 1} |\mu(f)| \le \sup_{\|f\|_{\tt L^k}\le 1} \sum_{z\in\Z^d} \mu(|f1_{B(z,r)}|)\le c(r) \sum_{z\in\Z^d} \|1_{B(z,r)}\ell_\mu\|_{L^{\ff k{k-1}}}.$$
 To prove the lower bound estimate in \eqref{SJ0},
  for each $z\in\Z^d$, we choose
  $f_z\in \B^+(\R^d)$ with $\|f_z1_{B(z,r)}\|_{L^k}= 1$ such that
  $$\mu(f_z  1_{B(z,r)})=\|\ell_{\mu}1_{B(z,r)}\|_{L^{\ff k{k-1}}}=\sup_{\|g\|_{L^k}\le 1}|\mu(g1_{B(z,r)})|.$$ This and
    \eqref{dfk} yield  that the function
  $$f:=\sum_{z\in\Z^d} f_z 1_{B(z,r)}$$ satisfies
$\|f\|_{\tt L^k}\leq c(r)$, so that
  $$c(r)^{-1} \sum_{z\in\Z^d} \|\ell_\mu 1_{B(z,r)}\|_{L^{\ff k{k-1}}}\le \|\mu\|_{k*}. $$
 Similarly, for any $\mu,\nu\in \scr P_{k*}$, we have
\beq\label{SJ}  c(r)^{-1} \sum_{z\in\mathbb Z^d} \|1_{B(z, r)} (\ell_\mu-\ell_\nu)\|_{L^{\ff k{k-1}}}
 \le \|\mu-\nu\|_{k*}\le c(r) \sum_{z\in\mathbb Z^d} \|1_{B(z, r)} (\ell_\mu-\ell_\nu)\|_{L^{\ff k{k-1}}}.  \end{equation}
From this  we see that $(\scr P_{k*},\ \|\cdot\|_{k*})$ is complete. Moreover,
since  $C_b(\R^d)$ is dense in $L^{k}(B(z,r))$ for any $z\in\Z^d$,
we may choose $\{f_n\}_{n\ge 1}\subset C_b(\R^d)$ such that
$$\|1_{B(z, r)} (\ell_\mu-\ell_\nu)\|_{L^{\ff k{k-1}}}= \sup_{n\ge 1}1_{\{\|f_n1_{B(z,r)}\|_{L^k}>0\}} \ff{|\mu(f_n)-\nu(f_n)|}{\|f_n1_{B(z,r)}\|_{L^k}},\ \ z\in\Z^d.$$
Combining this with \eqref{SJ}, we conclude that the Borel $\si$-field on $\scr P_{k*}$ induced by $\|\cdot\|_{k*}$ is contained by
that induced by the weak topology. Since the convergence in $\|\cdot\|_{k*}$ implies the weak convergence,
the former also contains the later, so that these two $\si$-fields  coincide each other.

(2) It suffices to prove for $\ll=0.$   Let $\{\mu^{(n)}\}_{n\ge 1}$ be
a Cauchy sequence with respect to $\rr_{T}^{p,k}$. Let $\oo_d$ be the volume of unit ball in $\R^d$. We have
$$\|f\|_{\tt L^k} \le \oo_d^{\ff 1 k} \|f\|_\infty,$$
so that $\|\cdot\|_{k*}\ge \oo_d^{-\ff 1 k} \|\cdot\|_{var} $ holds for the total variation norm $\|\cdot\|_{var}.$
By the completeness of $\|\cdot\|_{var}$ which is stronger than the weak topology,
there exists a unique $\mu\in C^w([0,T];\scr P)$ such that $\mu_0=\gg$ and
$$\lim_{n\to\infty} \|\mu_t^{(n)}-\mu_t\|_{var} =0,\ \ t\in [0,T].$$
Hence, for any $f\in \B_b(\R^d)$,
$$ |(\mu_t^{(n)}-\mu_t)(f)|=\liminf_{m\to\infty} |(\mu_t^{(n)}-\mu_t^{(m)})(f)|
\le \|f\|_{\tt L^k}\liminf_{m\to\infty} \|\mu_t^{(n)}-\mu_t^{(m)}\|_{k*},\ \ t\in [0,T].$$
This implies $$\|\mu_t^{(n)}-\mu_t\|_{k*}\le \liminf_{m\to\infty} \|\mu_t^{(n)}-\mu_t^{(m)}\|_{k*},\ \ t\in [0,T],$$
so that
$$\lim_{n\to\infty} \sup_{t\in (0,T]} t^{\ff {d(p-k)} {2pk}} \|\mu_t^{(n)}-\mu_t\|_{k*}
\le \lim_{m,n\to\infty} \sup_{t\in (0,T]} t^{\ff {d(p-k)} {2pk}}\|\mu_t^{(n)}-\mu_t^{(m)}\|_{k*}=0.$$
\end{proof}

 \beg{lem}\label{L2} Assume $(A_1)$  and   let $b$ be in $\eqref{B}$. Then for any  $T\in (0,\infty)$ and  $\mu\in\C_{p,k}^{\gg,T}$, the SDE $\eqref{EM}$ is weakly well-posed.
  If $(A_2)$  holds,
 then $\eqref{EM}$ is strongly well-posed.
 \end{lem}

 \beg{proof}   By $(A_1)$ and $\mu\in \C_{p,k}^{\gg,T}$, there exists a constant $c\in (0,\infty)$ such that
 $b_t^{0,\mu}(x):=b_t^{(0)}(x,\mu_t)$ satisfies
 $$|b_t^{0,\mu}(x,\mu_t)|\le c t^{-\ff{d(p-k)}{2pk}},\ \ t\in (0,T].$$
 Since $(p,k)\in \D$ implies $\ff{d(p-k)}{pk}<1$, we find $(p',q')\in \scr K$ such that
 $\|b^{0,\mu}\|_{\tt L_{q'}^{p'}(T)}<\infty.$
Then the desired assertions follows from  Proposition \ref{P01}.
 \end{proof}

By Lemma \ref{L2}, to confirm item (ii) above, it remains to verify \eqref{PHG}.   To this end, we introduce  local hyperbound estimates on the diffusion semigroup
\beq\label{SM}\bar P_{s,t}^\mu f(x):= \E[ f(\bar X_{s,t}^{\mu,x})],\ \ 0\le s\le t\le T,\ f\in \B_b(\R^d), x\in\R^d\end{equation}
for $\mu\in \C_{p,k}^{\gg,T}$, where    $\bar{X}_{s,t}^{\mu,x}$ (weakly) solves the SDE
\beq\label{BX} \d \bar X_{s,t}^{\mu,x}= \big\{b_t(\bar X_{s,t}^{\mu,x},\mu_t)- b_t^{(0)}(\bar X_{s,t}^{\mu,x},\mu_t)\big\}\d t+ \si_t(\bar X_{s,t}^{\mu,x})\d W_t,\ \ t\in [s,T],\ \bar X_{s,s}^{\mu,x}=x.\end{equation}
The next lemma follows from Proposition \ref{P01'}.

\beg{lem}\label{L3'} Assume $(A_1)$  and   let $b$ be in $\eqref{B}$. Then for any $T\in (0,\infty)$ and $1<p_1\le p_2\le\infty$, there exists a constant $c\in (0,\infty)$ such that for any $\gg\in \scr P_{p*}$ and $ \mu\in \C_{p,k}^{\gg, T}$,
\beq\label{ES1'} \|\bar P_{s,t}^\mu\|_{\tt L^{p_1}\to\tt L^{p_2}}\le c(t-s)^{-\ff{d(p_2-p_1)}{2p_1p_2}},
\ \ 0\le s\le t\le T,\end{equation}
\beq\label{ES2'} \|\nn\bar P_{s,t}^\mu\|_{\tt L^{p_1}\to\tt L^{p_2}}\le c(t-s)^{-\ff 1 2-\ff{d(p_2-p_1)}{2p_1p_2}},
\ \ 0\le s\le t\le T.\end{equation} When $ b^{(i)}=0$ for $2\le i\le l_0$, these estimates also hold for $p_1=1.$  \end{lem}

We are now ready to characterize the map   $\Phi^\gg$ defined in $\eqref{PHG}$ for $T=\tau_n(\gg)$.

 \beg{lem}\label{L3} Assume $(A_1)$  and   let $b$ be in $\eqref{B}$.
 Then the following assertions hold.
 \beg{enumerate}
 \item[$(1)$] Assume $p=\infty$. For any $n\in  \mathbb N$, there exist    $\ll(n),  \beta_1(n)\in (0,\infty)$ such that for any $\gg\in \scr P$,
   we have
\beq\label{MP}\Phi^{\gg}:\hat{\C}_{p,k}^{\gg, n,\lambda(n)}\to \hat{\C}_{p,k}^{\gg, n,\lambda(n)},\end{equation} where $\Phi^\gg$ is defined in $\eqref{PHG}$ for $T=n$ and
 \beq\label{TC0} \beg{split}&\hat{\C}_{p,k}^{\gg, n,\lambda(n)}:=\bigg\{\mu\in \C_{p,k}^{\gg, n}:\  \sup_{t\in (0, n]} \e^{-\ll(n)t} t^{\ff {d(p-k)}{2pk} }\|\mu_t\|_{k*}\le \beta_1(n) \bigg\}.\end{split}\end{equation}
For any $\Phi^\gg$-fixed point  $\mu\in \C_{p,k}^{\gg, n},$ we have
   $\mu\in \hat{\C}_{p,k}^{\gg, n,\lambda(n)}.$
\item[$(2)$] For any $n\in  \mathbb N$, there exist   constants $\beta_0(n)\in (0,1]$ and $  \beta_1(n)\in (0,\infty)$ such that for any $\gg\in \scr P_{p*}$ and $ \tau_n(\gg)$ defined in $\eqref{TT0}$,
   we have
\beq\label{MP1}\Phi^{\gg}: \tt\C_{p,k}^{\gg, n}\to \tt\C_{p,k}^{\gg, n},\end{equation} where $\Phi^\gg$ is defined in $\eqref{PHG}$ for $T=\tau_n(\gg)$ and
 \beq\label{TC01} \beg{split}&\tt\C_{p,k}^{\gg, n}:=\bigg\{\mu\in \C_{p,k}^{\gg, \tau_n(\gg)}:\  \sup_{t\in (0, \tau_n(\gg)]} t^{\ff {d(p-k)}{2pk} }\|\mu_t\|_{k*}\le \beta_1(n) \|\gg\|_{p*}\bigg\}.\end{split}\end{equation}
 For any $\Phi^\gg$-fixed point  $\mu\in \C_{p,k}^{\gg, \tau_n(\gg)},$ we have
   $\mu\in \tt\C_{p,k}^{\gg, n}.$
    \end{enumerate}
   \end{lem}

 \beg{proof}    We first prove that for fixed $T\in (0,\infty)$,
 \beq\label{PM} \Phi^\gg: \C_{p,k}^{\gg,T}\to \C_{p,k}^{\gg,T}.\end{equation} All constants $\{c_i:i\ge 0\}\subset (0,\infty)$ below do not depend on $\mu\in \C_{p,k}^{\gg,T}$.

 For $\mu\in \C_{p,k}^{\gg,T}$, let $\bar X_{s,t}^{\mu,x}$ solve \eqref{BX}, and denote $\bar X_t^{\mu,x}=\bar X_{0,t}^{\mu,x}$. Moreover, let    $X_t^{\mu,x}$ solve \eqref{EM} for $X_0^{\mu,x}=x,$  and let
\beq\label{PM'} P_t^\mu f(x):=\E[f(X_t^{\mu,x})],\ \ t\in [0,T],\ x\in\R^d,\ f\in \B_b(\R^d).
\end{equation}  By the definitions of $\|\cdot\|_{k*}$ and $\Phi^\gg\mu$, we have
\beq\label{K*} \|\Phi^\gg_t\mu\|_{k*}= \sup_{\|f\|_{\tt L^k}\le 1} \big|\gg(P_t^\mu f)\big|,\ \  t\in (0,T].\end{equation}
Noting that $(A_1)$ and $\mu\in \C_{p,k}^{\gg,T}$  imply that
$\xi_s:= (\si_s^*a_s^{-1} b_s^{(0)}) (\bar X_s^{\mu,x},\mu_s)$ satisfies
$$|\xi_s| \le  c_0 \rr_{T}^{p,k}(\mu)s^{-\ff {d(p-k)}{2pk}},\ \ s\in (0,T]$$ for some constant $c_0\in (0,\infty)$,
we see that
$$R_t:=\e^{\int_0^{t} \<\xi_s,\d W_s\>-\ff 1 2 \int_0^t|\xi_s|^2\d s},\ \ t\in [0,T]$$
is a martingale due to $(p,k)\in \D$. Noting that $k':=\ss k>1$, by  Girsanov's theorem, \eqref{K*} and \eqref{ES1'}, we find
a constant $C(\mu)\in (0,\infty)$ depending on $\mu$ such that
 \beg{align}\label{pgy}
 \nonumber&\|\Phi_t^\gg\mu\|_{k*}\leq\sup_{\|f\|_{\tt L^k}\le 1}\int_{\R^d} |\E [R_{t} f(\bar X_t^{\mu,x})]|\gamma(\d x)\\
&\le \|\gamma\|_{p*}\sup_{\|f\|_{\tt L^k}\le 1}\left\| \big(\E [R_{t}^{\ff {k'}{k'-1}}]\big)^{\ff {k'-1} {k'}}\big(\E[ |f|^{k'}(\bar X_t^{\mu,x})]\big)^{\ff 1 {k'}}\right\|_{\tilde{L}^p}\\
\nonumber&\le C(\mu) \|\gg\|_{p*} \|\bar P_t^\mu\|_{\tt L^{k'}\to \tt L^{p/k'}}^{1/k'}
\le C(\mu)c t^{-\ff{d(p-k)}{2pk}}\|\gamma\|_{p*},\ \ t\in (0,T].\end{align}
Hence, \eqref{PM} holds.

   From now on, let $T=\tau_n(\gamma)$  be in $\eqref{TT0}$ for some constant $\bb_0(n)\in (0,1]$ to be determined. By the Duhamel formula, see Proposition \ref{L*2}(2),
\begin{align}\label{DH1}  P_{r,t}^\mu f= \bar P_{r,t}^\mu f+\int_r^t P_{r,s}^\mu \<b_s^{(0)}(\cdot,\mu_s), \nn \bar P_{s,t}^\mu f\>\d s, \ \ 0\leq r\leq t\leq T,
\end{align}
we obtain 
\beq\label{DH}\beg{split}&(\Phi_t^\gg \mu)(f)=\gg(P_t^\mu f)   \\
&= \gg(\bar P_{t}^\mu  f) + \int_0^t \gg\big(P_s^\mu \<b_s^{(0)}(\cdot,\mu_s), \nn   \bar P_{s,t}^\mu  f\> \big)\d s,\ \ t\in [0, T].\end{split}\end{equation}
 Below we consider three different  cases respectively:  1) $p=\infty$;  2) $b^{(0)}=0$, and  3)  $p<\infty$ with $ b^{(0)}\ne 0$. All constants below may depend on $n$.

Having the above preparations, we are able to prove assertions (1) and (2) in three different cases.

(1) Assume $p=\infty$. In this case,  $ T:=\tau_n(\gg)=n$.
   Since $\|P_t^\mu\|_{\tt L^\infty\to\tt L^\infty}=1,$ by $(A_1)$ for $T=n$,   \eqref{ES1'}, \eqref{ES2'} and \eqref{DH}, we find  a constant  $c_1\in (0,\infty)$ such that
  \beg{align*} &\|\Phi_t^\gg\mu\|_{k*}= \sup_{\|f\|_{\tt L^k}\le 1} |(\Phi_t^\gg \mu)(f)|\\
 &\le \|\bar P_t^\mu\|_{\tt L^k\to \tt L^\infty}+ K\int_0^t \|P_s^\mu\|_{\tt L^\infty\to\tt L^\infty}\|\mu_s\|_{k*}\|\nn \bar P_{s,t}^\mu\|_{\tt L^k\to\tt L^\infty}\d s\\
 &\le c_1 t^{-\ff {d}{2k}}  + c_1 \int_0^t \|\mu_s\|_{k*} (t-s)^{-\ff 1 2-\ff d{2k}}\d s.\end{align*} So, there exist  constants $c_2,c_3\ge 1$ such that for any $\ll\in (0,\infty)$,
 \beq\label{U}\beg{split} &\rr_{n,\ll}^{p,k}(\Phi^\gg\mu):=\sup_{t\in [0,n]}t^{\ff{d}{2k}}\e^{-\ll t} \|\Phi_t^\gg\mu\|_{k*}\\
 &\le c_2 + c_1 \rr_{n,\ll}^{p,k}(\mu)  \sup_{t\in (0,n]}t^{\ff d{2k}}\int_0^t s^{-\ff {  d}{2k}}\e^{-\ll(t-s)}(t-s)^{-\ff 1 2-\ff d{2k}}\d s
 \\
 &\le c_3    + c_3 \rr_{n,\ll}^{p,k}(\mu) \ll^{-\theta_0},\ \ \theta_0:= \ff 1 2-\ff { d}{2k}>0,\ \ t\in (0,n].\end{split}\end{equation}
Letting
\beq\label{LL} \ll(n):= (2c_3)^{\theta_0^{-1}}, \end{equation} we obtain
$$\rr_{n,\ll}^{p,k}(\Phi^\gg\mu)\le 2c_3,\ \ \text{if}\ \ \rr_{n,\ll}^{p,k}(\mu)\le 2 c_3.$$ Noting that
$\|\mu\|_{\infty*}=1$ for $\mu\in \scr P$,  we conclude that  \eqref{MP}   holds  for
 $\beta_1(n):=   c_3,  \ll(n)=(2c_3)^{\theta_0^{-1}}, \tau_n(\gg)=n $   and $\tt\C_{\infty,k}^{\gg,n}$   in \eqref{TC0} with $p=\infty$.

If $\mu$ is a fixed point of $\Phi^\gg$ such that
 $\Phi_t^\gg\mu=\mu_t$,  then \eqref{U} and \eqref{LL} imply
  $$\rr_{n,\ll(n)}^{p,k}(\mu)\le c_3    + c_3\rr_{n,\ll(n)}^{p,k}(\mu)\ll^{-\theta_0}= c_3    + \ff 1 2 \rr_{n,\ll(n)}^{p,k}(\mu).$$
  So, $\mu\in \tt \C_{p,k}^{\gg,n}.$

(2) (i) Assume   $b^{(0)}=0$.    In this case, $ T:=\tau_n(\gg)=n$ and $P^\mu_t=\bar P_t^\mu$.  By \eqref{ES1'}, we find a constant $\bb_1(n)\in (0,\infty)$ such that
\begin{align*}
&\|\Phi_t^\gg\mu\|_{k*}= \sup_{\|f\|_{\tt L^k}\le 1} |(\Phi_t^\gg \mu)(f)|\leq\sup_{\|f\|_{\tt L^k}\le 1} \gamma(|\bar P_t^\mu f|)\\
 &\le \|\gamma\|_{p*}\sup_{\|f\|_{\tt L^k}\leq 1}\|\bar P_t^\mu f\|_{\tilde{L}^p}=\|\gamma\|_{p*}\|\bar P_t^\mu\|_{\tt L^k\to \tt L^p}\\
 &\leq \bb_1(n) \|\gamma\|_{p*}t^{-\frac{d(p-k)}{2pk}},\ \  \ t\in(0,n],\ \mu\in \C_{p,k}^{\gg,n}.
\end{align*}
Thus, \eqref{MP1} holds,  and any fixed point of $\Phi^\gg$ belongs to $\tt \C_{p,k}^{\gg,n}$ defined in \eqref{TC01}.

(ii)   $p<\infty$ and $ b^{(0)}\ne  0$.   By $(A_1)$ for $T=\tau_n(\gg)$, \eqref{ES1'}, \eqref{ES2'}, \eqref{K*} and \eqref{DH},  we find constants $c_1,c_2\in (0,\infty)$ such that
\beq\label{MTY} \beg{split} \|\Phi_t^\gg\mu\|_{k*}&\le c_1 t^{-\ff {d(p-k)}{2pk}}\|\gamma\|_{p*}+ K \int_0^t \|\Phi_s^\gg\mu\|_{k*}\|\mu_s\|_{k*}
 \sup_{\|f\|_{\tt L^k}\leq 1}\|\nn \bar P_{s,t}^\mu f\|_{\tt L^k}\d s\\
&\le  c_1 t^{-\ff {d(p-k)}{2pk}}\|\gamma\|_{p*}+ c_2  \int_0^t \|\Phi_s^\gg\mu\|_{k*}\|\mu_s\|_{k*}(t-s)^{-\ff 1 2} \d s,\ \ t\in (0,T].\end{split}\end{equation}
Noting that $\rr_t^{p,k}(\mu)$ is non-decreasing in $t$,
\beq\label{RL} \rr_{t+}^{p,k}(\mu):= \lim_{\vv\downarrow 0} \rr_{(t+\vv)\land T}^{p,k}(\mu),\ \
\rr_{t-}^{p,k}(\mu):= \lim_{\vv\downarrow 0} \rr_{(t-\vv)\vee 0}^{p,k}(\mu)\end{equation}
exist and are non-decreasing for $t\in (0,T].$
By \eqref{pgy} and $\mu\in \C_{p,k}^{\gg,T}$,
we find a constant $C'(\mu) \in (0,\infty)$ such that
\beq\label{AQ} \|\Phi_s^\gg\mu\|_{k*}\|\mu_s\|_{k*}\le \rr_s^{p,k}(\Phi^\gg\mu) \rr_s^{p,k}(\mu)s^{-\ff{d(p-k)}{pk}}\le  C'(\mu) s^{-\ff{d(p-k)}{pk}},\ \ s\in (0,T].\end{equation}
Since    $\ff{d(p-k)}{pk}<1$ due to $(p,k)\in\D$,  \eqref{AQ} implies that the function
$$(0,T]\ni t\mapsto \int_0^t \|\Phi_s^\gg\mu\|_{k*}\|\mu_s\|_{k*}(t-s)^{-\ff 1 2} \d s$$ is continuous. Combining this with \eqref{MTY}, \eqref{RL} and \eqref{AQ}, we  find  constants $c_3, c_4\ge 1$ such that
\beq\label{I}\beg{split}  &\rr_{t+}^{p,k}(\Phi^\gg\mu):= \lim_{\vv\downarrow 0}\sup_{s\in (0,(t+\vv)\land T]} s^{\ff {d(p-k)}{2pk}} \|\Phi_s^\gg\mu\|_{k*}\\
&\le c_1\|\gg\|_{p*} + c_2\sup_{s\in (0,t]}s^{\ff {d(p-k)}{2pk}} \int_0^s \|\Phi_r^\gg\mu\|_{k*}\|\mu_r\|_{k*}(s-r)^{-\ff 1 2}\d r\\
&\le c_1 \|\gamma\|_{p*}+c_3 \rr_{t-}^{p,k}(\Phi^\gg\mu) \rr_{t-}^{p,k}(\mu) \sup_{s\in (0,t]} s^{\ff {d(p-k)}{2pk}}\int_0^sr^{-\ff {d(p-k)}{pk}}(s-r)^{-\ff 1 2}\d r\\
&\le c_4 \|\gg\|_{p*} +c_4  \rr_{t-}^{p,k}(\Phi^\gg\mu) \rr_{t-}^{p,k}(\mu) t^{\theta},\ \  t\in (0, T], \end{split}\end{equation} where $\theta:= \ff 1 2 - \ff{d(p-k)}{2pk}>0. $
Letting $\bb_1(n)=2c_4$, we obtain
\beq\label{I'} \beg{split} \rr_{t+}^{p,k}(\Phi^\gg\mu)\le c_4  \|\gg\|_{p*} + 2c_4^2   \rr_{t-}^{p,k}(\Phi^\gg\mu)  \|\gg\|_{p*} t^{\theta},\
   \  t\in  (0, T],\ \ \mu\in \tt\C_{p,k}^{\gg, n}.\end{split} \end{equation}
 So,  for $ T=\tau_n(\gg)$ in \eqref{TT0} with
$  \beta_0(n):= (4c_4^2)^{-1/\theta},$  we have
$$ 2c_4^2  \|\gg\|_{p*}  t^{\theta}\le  2c_4^2\beta_0(n)^{\theta}  = \ff 1 2,\ \ \ t\in (0, T].$$
Hence, \eqref{I'} implies \eqref{MP1}.

If $\mu\in \C_{p,k}^{\gg, T}$ is a fixed point of $\Phi^\gg$, then $\Phi^\gg\mu=\mu$ so that
\eqref{I} implies
\beq\label{PI}\beg{split}\rr_{t+}^{p,k}(\mu)\le c_4   \|\gg\|_{p*} + c_4 \rr_{t-}^{p,k}(\mu)^2 t^{\theta},\ \ t\in (0, T].\end{split}\end{equation}
Then
\beq\label{PKb} \rr_{0+}^{p,k}(\mu):=  \lim_{t\downarrow 0}  \rr_{t\wedge T}^{p,k}(\mu) \le c_4 \|\gg\|_{p*}.\end{equation}
 This and the right continuity of $\rr^{p,k}_{t+}$ in $t$ imply
$$s_0:=  T\land \inf\big\{t\in (0, T]:\  \rr_{t+}^{p,k}(\mu)\ge 2 c_4 \|\gg\|_{p*}\big\}>0,$$
where $\inf\emptyset:=\infty$ by convention.
 If $s_0< T$,  by the non-decreasing of $\rr_{t}^{p,k}$ and \eqref{RL}, we obtain
 $$  \rr_{s_0+}^{p,k}(\mu)\ge 2 c_4 \|\gg\|_{p*}\ge \rr_{s_0-}^{p,k}(\mu),$$
 so that  \eqref{PI} yields
 \beg{align*}& 2 c_4 \|\gg\|_{p*} \le  \rr_{s_0+}^{p,k}(\mu)  \leq c_4 \|\gg\|_{p*}  +  4 c_4^3 \|\gg\|_{p*}^2   s_0^{\theta},\end{align*}
and thus,
$$ s_0\ge \big(4c_4^2\|\gg\|_{p*}  \big)^{-1/\theta}=  T,$$
which contradicts  to $s_0< T$. Hence, $s_0= T$, so that     \eqref{PI} together with $\rr_{s_0-}^{p,k}(\mu)\le  2 c_4 \|\gg\|_{p*}$ and
$$s_0^\theta=T^\theta=\tau_n(\gg)^\theta=\beta_0(n)^{\theta}\|\gg\|_{p*}^{-1}=(4c_4^2)^{-1}\|\gg\|_{p*}^{-1}$$
 implies
$$\rr_{ T}^{p,k}(\mu) \le c_4   \|\gg\|_{p*} + c_4 \rr_{s_0-}^{p,k}(\mu)^2 s_0^{\theta}= 2 c_4 \|\gg\|_{p*} =\beta_1(n)  \|\gg\|_{p*}.$$
Therefore, $\mu\in\tt\C_{p,k}^{\gg, n}.$
 \end{proof}

We are now ready to solve \eqref{E0} with initial value $\gg\in \scr P_{p*}$ up to  time $\tau_n(\gg)$ for any $n\in\mathbb N$.

\beg{prp}\label{T01} Assume $(A_1)$. Let $b$ be in $\eqref{B}$, and let $n\in\mathbb N$.
Then the following assertions hold.
\beg{enumerate} \item[$(1)$]  There exist constants
$\bb_0(n)\in (0,1]$ and $\bb_1(n)\in (0,\infty)$ such that for any $\gg\in \scr P_{p*}$, the SDE $\eqref{E0}$ with initial distribution $\gg$ has a unique weak $\C_{p,k}$-solution   up to time $\tau_n(\gg)$ defined
in $\eqref{TT0}$,
and $\eqref{EST}$ holds.
 \item[$(2)$]  If   $(A_2)$ holds,    then for any $\F_0$-measurable initial value $X_0$ with $\gg:=\L_{X_0}\in \scr P_{p*}$, the SDE $\eqref{E0}$ has a unique  strong $\C_{p,k}$-solution up to time $\tau_n(\gg)$, and for any $q\in [1,\infty)$  there exists a constant $c(n,q)\in (0,\infty)$,   such that  for any $\L_{X_0}=\gg\in \scr P_{p*}$,
\beq\label{NESN}  \E\bigg[\sup_{s\in [0,\tau_n(\gg)]} |X_s|^q\bigg|\F_0\bigg]\le c(n,q) (1+|X_0|^q).\end{equation}   \end{enumerate}  \end{prp}

\beg{proof}  Simply denote $ \tau_n=\tau_n(\gg)$ an let $\tt\C_{p,k}^{\gg,n}$ be in \eqref{TC01} .The following procedure is also available in the case $p=\infty$ by replacing $ \tau_n(\gg)$ and $\tt\C_{p,k}^{\gg,n}$ with $n$ and $\hat{\C}_{p,k}^{\gg, n,\lambda(n)}$ in \eqref{TC0} respectively. So, we only consider the case $p<\infty$.

(1)  By Lemma \ref{L3},
  all fixed points in $\C_{p,k}^{\gg, \tau_n}$ of $\Phi^\gg$ are included in $\tt\C_{p,k}^{\gg, n}$.
Therefore, \eqref{EST} holds for any (weak) $\C_{p,k}$-solutions of $\eqref{E0}$ with initial distribution $\gg$ up to time $ \tau_n$.
By Lemma \ref{L1} and the contractive  fixed point theorem, it suffices to find $\ll\in (0,\infty)$ such  that
$$\Phi^\gg:\ \tt\C_{p,k}^{\gg, n} \to \tt\C_{p,k}^{\gg, n}$$ is contractive under the metric  $\rr_{ \tau_n,\ll}^{p,k}$.

Let $\mu,\nu\in \tt\C_{p,k}^{\gg, n}.$ Since $\gg\in \scr P_{p*}$ is given, positive constants  in the following are  allowed  to depend on $\gg$.
Let $\L_{\bar X_0}=\gg$ and $\bar  X_t$ (weakly) solve the SDE
$$\d \bar X_t= (b_t-b_t^{(0)}) (\bar X_t, \mu_t)\d t+\si_t(\bar X_t)\d W_t,\ \ t\in [0, \tau_n].$$
Then
$$\E[f(\bar X_t)]= \gg(\bar P_t^\mu f),\ \ t\in [0, \tau_n], f\in \B_b(\R^d),$$
where $\bar P_t^\mu=\bar P_{0,t}^\mu $  is in \eqref{SM}.
Let
\beg{align*}&\xi_s^1:= (\si_s^*a_s^{-1})(\bar X_s) b_s^{(0)}(\bar X_s,\mu_s),\\
&\xi_s^2:= (\si_s^*a_s^{-1})(\bar X_s) \{b_s(\bar X_s,\nu_s)-b_s(\bar X_s,\mu_s)+b_s^{(0)}(\bar X_s,\mu_s)\},\ \ s\in [0, \tau_n].\end{align*}
By $(A_1)$ and $\mu,\nu\in \tt \C_{p,k}^{\gg,n}$, we find a constant $c_1\in (0,\infty)$  such that
\beq\label{C1}\beg{split}&|\xi_s^i|^2  \le c_1s^{-\ff {d(p-k)}{pk}},\ \ i=1,2,\\
&|\xi_s^1-\xi_s^2|^2 \le c_1 \|\mu_s-\nu_s\|_{k*}^2,\ \ \ \ s\in (0, \tau_n].\end{split}\end{equation}
Since $(p,k)\in \D$ implies $\ff {d(p-k)}{pk}<1$, and noting that  $\mu,\nu\in \tt\C_{p,k}^{\gg, n}$ implies
$$ \|\mu_s-\nu_s\|_{k*}^2\le c s^{-\ff {d(p-k)}{pk}}$$
for some constant $c \in (0,\infty)$,
by Girsanov's theorem,
$$R_t^i:=\e^{\int_0^t \<\xi_s^i,\d W_s\>-\ff 1 2 \int_0^t |\xi_s^i|^2\d s},\ \ t\in [0, \tau_n],\ i=1,2$$
are  martingales, and
\beg{align*} &\|\Phi_t^\gg\mu-\Phi_t^\gg\nu\|_{k*}=\sup_{\|f\|_{\tt L^k\le 1}} |\E[(R_t^1-R_t^2)f(\bar X_t)]|\\
&\le \sup_{\|f\|_{\tt L^k\le 1}} \E\bigg[\big(\E[|R_t^1-R^2_t|^{\ff {k}{k-1}}|\F_0]\big)^{\ff {k-1}{k}} \big(\E[|f|^{k}(\bar X_t)|\F_0]\big)^{\ff 1 {k}}\bigg].\end{align*}
By \eqref{C1}, we find   constants $c_2,c_3\in (0,\infty)$ such that
$$\big(\E[|R_t^1-R_t^2|^{\ff {k}{k-1}}|\F_0]\big)^{\ff {k-1}{k}}\le c_2  \left(\int_0^t\|\mu_s-\nu_s\|_{k*}^2\d s\right)^{\frac{1}{2}},$$
and by \eqref{ES1'},
\beg{align*}&\E\bigg[\big(\E[|f|^{k}(\bar X_t)|\F_0]\big)^{\ff 1 {k}}\bigg]=\gg\big((\bar P_t^\mu|f|^{k})^{\ff 1 {k}}\big)\\
&\le \|\gg\|_{p*} \|f\|_{\tt L^k} \|\bar P_t^\mu\|_{\tt L^{1}\to \tt L^{\ff p k}}^{\ff 1 k}
\le c_3   \|f\|_{\tt L^k}t^{-\ff{d(p-k)}{2pk}},\ \ t\in (0, \tau_n].\end{align*}
Therefore, there exists a constant $c_4\in (0,\infty)$ such that for any $\ll\in (0,\infty)$,
 \beg{align*}\rr_{ \tau_n,\ll}^{p,k}(\Phi^\gg\mu,\Phi^\gg\nu)&\le c_2 c_3
 \rr_{\tau_n,\ll}^{p,k}(\mu,\nu)\sup_{t\in (0, \tau_n]}\left(\int_0^{t}  s^{-\ff{d(p-k)}{pk}}\e^{-2\ll(t-s)}\d s\right)^\frac{1}{2}\\
&\le c_4 \ll^{\ff{d(p-k)}{2pk}-\frac{1}{2}} \rr_{ \tau_n,\ll}^{p,k}(\mu, \nu).\end{align*}
Since $\ff{d(p-k)}{pk}-1<0$ due to $(p,k)\in \D$, when $\ll\in (0,\infty)$ is large enough $\Phi^\gg$ is contractive under
$\rr_{ \tau_n,\ll}^{p,k}$ as desired.

(2) By Lemma \ref{L2}, if  {\bf (A)} holds,   then $\eqref{EM}$ for $T=\tau_n$ is strongly well-posed for any $\mu\in \C_{p,k}^{\tau_n}$. Combining this with the weak well-posedness of \eqref{E0} ensured by  Proposition \ref{T01}(1),
we derive the strong well-posedness of \eqref{E0} up to time $\tau_n$.

To prove \eqref{NESN}, we   consider the SDE
$$\d \bar X_t= \big\{b_t(\bar X_t, P_t^*\gg)- b^{(0)}_t(\bar X_t, P_t^*\gg)\big\}\d t+\si_t(\bar X_t)\d W_t,\ \ \bar X_0=X_0,\ t\in [0, \tau_n].$$
According to Proposition \ref{L*1},  {\bf (A)}  implies that this SDE is well-posed and there exists a constant $c_1(n,q)\in (0,\infty)$  independent of the initial distribution $\gg$
such that
\beq\label{NE1} \E\bigg[\sup_{t\in [0, \tau_n]} |\bar X_t|^q\bigg|\F_0\bigg]\le c_1(n,q) (1+|X_0|^q).\end{equation}
When $b^{(0)}=0,$ we have $\tau_n(\gg)=n$ and $X_t=\bar X_t,$ so that \eqref{NESN} holds.

For $b^{(0)}\ne 0$, let
$$\xi_t:=(\si_t^*a_t^{-1})(\bar X_t)  b_t^{(0)} (\bar X_t, P_t^*\gg),\ \ \ t\in [0, \tau_n].$$
By $(A_1)$, \eqref{TT0}, $ \tau_n=\tau_n(\gg)$  and \eqref{EST}, we find   constants $k_1,k_2  \in (0,\infty)$ such that
 \beg{align*} \int_0^{ \tau_n} |\xi_t|^2 \d t& \le k_1+k_1 (1+\|\gg\|_{p*})^2 \int_0^{ \tau_n}    t^{-\ff {d(p-k)}{pk}} \d t\\
&\le  k_1+k_2 (1+\|\gg\|_{p*})^2  \tau_n^{1-\ff {d(p-k)}{pk}} = k_1 +k_2\bb_0(n)^{-2}=:k_3.  \end{align*}
 So,
$$R_t:=\e^{\int_0^t\<\xi_s,\d W_s\>- \int_0^t \frac{1}{2}|\xi_s|^2\d s},\ \ \ t\in [0, \tau_n]$$
is an exponential martingale, and by  Girsanov's theorem and \eqref{NE1},  we obtain
\beg{align*}  &\E\bigg[\sup_{t\in [0, \tau_n]} |X_t|^q\Big|\F_0\bigg]=  \E\bigg[R_{ \tau_n} \sup_{t\in [0, \tau_n]} |\bar X_t|^q\Big|\F_0\bigg]\\
&\le \big(\E[R_{ \tau_n}^2|\F_0]\big)^{\ff 1 2}  \bigg(\E\bigg[\sup_{t\in [0, \tau_n]} |\bar X_t|^{2q}\Big|\F_0\bigg]\bigg)^{\ff 1 2}
 \le \e^{k_3}\ss{c_1(n,2q) (1+|X_0|^{2q})}.   \end{align*}
Therefore, \eqref{NESN} holds for some constant $c(n,q)\in (0,\infty).$

  \end{proof}

\beg{proof}[\textbf{Proof of Theorem \ref{T0}}] Let $\gg\in \scr P_{p*}.$

(a) If $\tau_n(\gg)=n$ holds for any $n\in \mathbb N$,  then by Proposition \ref{T01},  the SDE \eqref{E0} with initial distribution $\gg$ has a unique weak $\C_{p,k}$-solution up to any time $t\ge 0$, so that $\tau(\gg)=\infty$.

If there exists $n\in \mathbb N$ such that $\tau_n (\gg)<n$,
  by applying Proposition \ref{T01} to the SDE \eqref{E0} starting from time $\tau_n(\gg)$   with initial distribution $\gg_0:=\L_{X_{\tau_n(\gg)}}$,
we conclude that \eqref{E0} has a unique  weak  $\C_{p,k}$-solution up to time
$$\tau_{n,1}(\gg):=n\land \big(\tau_n(\gg)+   \beta_0 (n) \|\gg_0\|_{p*}^{ -\ff{1}{\theta}}\big). $$  In general, once   \eqref{E0} has a unique  weak $\C_{p,k}$-solution up to time
$ \tau_{n,i} (\gg)$ for some $ i\in \mathbb N$ so that $\gg_i:=\L_{X_{\tau_{n,i}(\gg)}}\in \scr P_{k*}$,  it also  has a unique  weak  $\C_{p,k}$-solution up to time
$$\tau_{n,i+1}:=n\land \big(\tau_{n,i}(\gg)+   \beta_0 (n) \|\gg_{i}\|_{p*}^{ -\ff{1}{\theta}}\big).$$
 Hence,      we find a deterministic life time
$$\hat\tau_n(\gg):=\lim_{i\to\infty} \tau_{n,i} \in (\tau_n(\gg),n]$$ such that
\eqref{E0} has a unique  weak $\C_{p,k}$-solution up to  any  time $t\in [\tau_n(\gg),\hat\tau_n(\gg)),$
and when $\hat\tau_n(\gg)<n$
  $$\limsup_{t\to \hat\tau_n(\gg)} \|\L_{X_t}\|_{p*}=\infty.$$
Let
  $ \tau(\gg)=\hat\tau_n(\gg)$ for the smallest $n\in\mathbb N$ with $\hat\tau_n(\gg)<n$, and let $\tau(\gg)=\infty$ if such $n$ does not exist. Then
 \eqref{E0} has a unique maxial weak $\C_{p,k}$-solution with life time $\tau(\gg)$.
We have proved  Theorem \ref{T0}(1)-(2)
since \eqref{EST} is included by Proposition \ref{T01}.

(b) If $\tau(\gg)<\infty$, then $\tau(\gg)<n$ for some  $n\in \mathbb N$.    If \eqref{L} does not hold, then for any $\vv\in (0,1)$ we find $(0, \tau(\gg)) \ni \vv_i\downarrow 0$ as $i\uparrow\infty$  such that for $s_i:= \tau(\gg)-\vv_i$ satisfies
$$\|P_{s_i}^*\gg\|_{p*}\le \vv \vv_i^{-\theta},\ \ \ i\ge 1.$$
By  Proposition \ref{T01}  for $\eqref{E0}$ starting from time $s_i\le n$, we conclude that this SDE has a unique weak $\C_{p,k}$-solution up to time
$$ s_i+  \bb_0 (n)(\vv \vv_i^{-\theta})^{-\theta^{-1}},\ \ i\ge 1.$$
So,
$$\tau(\gg) \ge   s_i+  \bb_0 (n)( \vv \vv_i^{-\theta})^{-\theta^{-1}}= \tau(\gg)-\vv_i +  \bb_0 (n) (\vv \vv_i^{-\theta})^{-\theta^{-1}},\ \ i\ge 1.$$
Thus,
$$1\le \lim_{i\to\infty} \vv_i  \bb_0 (n)^{-1}  (\vv \vv_i^{-\theta})^{\theta^{-1}} =  \bb_0(n)^{-1} \vv^{\theta^{-1}},$$
which contracts to the arbitrariness of $\vv\in (0,1).$ Hence \eqref{L} holds.

Next, let $\mu_t:=P_t^*\gg, t\in [0,\tau(\gg)),$  let $\bar X_t^\mu$ solve the SDE \eqref{BX} for $\L_{\bar X_0^\mu}=\gg$, and let $\bar P_t^\mu=\bar P_{0,t}^\mu$ be in \eqref{SM}. Then
\beq\label{99} \E[f(\bar X_t^\mu)]= \gg(\bar P_t^\mu f), \ \
t\in (0,\tau(\gg)),\|f\|_{\tt L^p}<\infty.\end{equation}
By Girsanov's theorem, we have
\beq\label{90} (P_t^*\gg)(f)= \E[f(\bar X_t^\mu)R_t],\ \  t\in (0,\tau(\gg)),\|f\|_{\tt L^p}<\infty,\end{equation}
where  $R_t:=\exp[\int_0^t\<\zeta_s,\d W_s\>-\ff 1 2\int_0^t |\zeta_s|^2\d s]$ for
$$\zeta_s:=  \big(\si_s^*a_s^{-1}\big)(\bar X_s^\mu) b_s^{(0)}(\bar X_s^\mu,\mu_s).$$
By $(A_1)$, we find a constant $K\in (0,\infty)$ such that
$|\zeta_s|\le K \|\mu_s\|_{k*}$ for $s\in (0,\tau(\gg))$. Hence, for any  $\aa\in (1,p)$, we find  a constant  $c_1\in (0,\infty)$ such that
\beq\label{92} \E[R_t^{\ff\aa{\aa-1}}] \le \e^{c_1\int_0^t \|\mu_s\|_{k*}^2\d s}.\end{equation}
Combining \eqref{99}-\eqref{92} and H\"older's inequality, we derive
\beg{align*} &\big| (P_t^*\gg)(f)\big|\le \e^{c_1\int_0^t \|\mu_s\|_{k*}^2\d s} \big(\E[|f|^\aa (\bar X_t^\mu)]\big)^{\ff 1 \aa}\\
& =  \e^{c_1\int_0^t \|\mu_s\|_{k*}^2\d s} \big[\gg(\bar P_t^{\mu} |f|^\aa)\big]^{\ff 1\aa} \le \e^{c_1\int_0^t \|\mu_s\|_{k*}^2\d s}  \|f\|_{\tt L^p}\big( \|\gg\|_{p*}  \|\bar P_t^\mu\|_{\tt L^{\ff p\aa}\to \tt L^p}\big)^{\ff 1 \aa}.\end{align*}
This together with \eqref{ES1'} implies that
$$\|P_t^*\gg\|_{p*}= \sup_{\|f\|_{\tt L^p}\le 1} \big| (P_t^*\gg)(f)\big|\le   \e^{c_1\int_0^t \|\mu_s\|_{k*}^2\d s}\Big(c\|\gg\|_{p*} t^{-\ff{d(\aa-1)}{2p}}\Big)^{\ff 1 \aa},\ \ t\in (0,\tau(\gg)).$$
Since  $\limsup_{t\uparrow\tau(\gamma)} \|P_t^*\gg\|_{p*}=\infty$  due to  \eqref{L},  we obtain
$$ \int_0^{\tau(\gg)} \|\mu_s\|_{k*}^2\d s  =\infty.$$
 Therefore,  \eqref{L'} holds for any   $r\in (0,\tau(\gg)),$ since by the definition of maximal $\C_{p,k}$-solution we find a constant $c(r)\in (0,\infty)$ such that
$$\int_0^r \|\mu_s\|_{k*}^2   \d s \le   c(r) \|\gg\|_{p*}^2  \int_0^r s^{-\ff{d(p-k)}{pk}}\d s <\infty,$$
where $\ff{d(p-k)}{pk}<1$ by $(p,k)\in \D.$
\end{proof}

 \beg{proof}[\textbf{Proof of Corollary \ref{C01}}]  Let $\gg\in \scr P_{p*}$ and $s\in (0,\tau(\gg))$.

 (a) If \eqref{E0} has a weak $\C_{p,k}$-solution $(\tt X_t)_{t\in [0,T]}$ up to a finite time $T>s$, then $\tau(\gg)>T$.
Indeed, by the weak uniqueness up to time $\tau(\gg)$ due to   Theorem \ref{T0},
we have $\L_{X_t}=\L_{\tt X_t}$ for $t<T\land \tau(\gg)$, which  together with $s<T\land\tau(\gg)$ and $\L_{\tt X_\cdot}\in \C_{p,k}^T$  implies
$$\limsup_{t\uparrow T\land\tau(\gg) } \|\L_{X_t}\|_{k*}\le \sup_{t\in [s,T]} \|\L_{\tt X_t} \|_{k*}<\infty,$$
so that   $\tau(\gg)> T$ according to Theorem \ref{T0}(3).

(b) Denote  $\mu_t:= P_t^*\gg, \ t\in [0,\tau(\gg))$. We first prove $\|\mu_{t}\|_{1*}<\infty$ for any $t\in (0,\tau(\gg))$.
By $(A_1)$ and \eqref{EST}, we find a constant $c_1(t)\in (0,\infty)$ depending on $\gg$ and increasing in $t\in (0, \tau(\gg))$  such that
\beq\label{BNP} \|b^{(0)}_s(\cdot,\mu_s)\|_{\infty}\le c_1(t) s^{-\ff {d(p-k)}{2pk}},\ \ s\in (0,t].\end{equation}
So,   there exists $(p_0',q_0')\in \scr K$   such that $\|b^{(0)}_\cdot(\cdot,\mu_\cdot)\|_{\tt L_{q_0'}^{p_0'}(s,t)}<\infty.$ By Lemma \ref{L3'} for $l_0=2$ and $b_t^{(0)}(\cdot,\mu_t)$ in place of $b_t^{(2)}(\cdot,\mu_t)$, we derive \eqref{ES1'} and \eqref{ES2'} for
$P_t^\mu$ in place of   $\bar P_{s,t}^\mu.$ So,   for fixed  $l\in (1,p\land \ff{d}{(d-1)^+})$,
we find a constant $c(l,t)\in (0,\infty)$ increasing in $t$  such that
\beq\label{9} \beg{split} &\|\mu_t\|_{1*}=\sup_{\|f\|_{\tt L^1}\le 1}|\mu_t(f)|= \sup_{\|f\|_{\tt L^1}\le 1}|\gg (P_t^\mu f)|\le \|\gg\|_{p*} \|P_t^\mu\|_{\tt L^1\to\tt L^p}\\
&= \|\gg\|_{p*} \big\| P_{\ff t 2}^\mu  P_{\ff t 2,t}^\mu\big\|_{\tt L^1\to\tt L^p}
 \le  \|\gg\|_{p*} \| P_{\ff t 2}^\mu\|_{\tt L^l\to\tt L^p} \| P_{\ff t 2,t}^\mu\|_{\tt L^1\to\tt L^l} \\
 & \le c(l,t) t^{-\ff{d(p-l)}{2pl}}
 \big\|P_{\ff t 2,t}^\mu\big\|_{\tt L^1\to\tt L^l},\ \ t\in (0,\tau(\gg)).\end{split}\end{equation}
By Lemma \ref{L3'} for $b^{(i)}=0, 2\leq i\leq l_0$,  \eqref{BNP}  and Duhamel's formula \eqref{DH1} for $r=\ff t 2$, i.e.
$$ P_{\ff t 2,t}^\mu f= \hat{ P}_{\ff t2,t} f+\int_{\ff t 2}^t  P_{\ff t 2,s}^\mu\<b_s^{(0)}(\cdot,\mu_s),\nn \hat P_{s,t} f\>\d s,$$
we find   constants $c_2(t),c_3(t),c_4(t)\in (0,\infty)$ increasing in $t$ such that
\beg{align*} & \big\| P_{\ff t 2,t}^\mu\big\|_{\tt L^1\to\tt L^l}\le c_2(t) t^{-\ff{d(l-1)}{2l}} +c_2(t)   t^{-\ff {d(p-k)}{2pk} }\int_{\ff t 2}^t \big\| P_{\ff t 2,s}^\mu\big\|_{\tt L^l\to\tt L^l}  \|\nn\hat P_{s,t}\|_{\tt L^1\to \tt L^l}\d s\\
&\le c_2(t) t^{-\ff{d(l-1)}{2l}} +c_3 (t)  t^{-\ff {d(p-k)}{2pk}} \int_{\ff t 2}^t(t-s)^{-\ff 1 2-\ff{d(l-1)}{2l}}\d s \le c_4(t) t^{-\ff{d(l-1)}{2l}},\ \ t\in (0,\tau(\gg)),\end{align*}
where the last step follows from $\ff {d(p-k)}{2pk}\le \ff 1 2$ and $\ff{d(l-1)}{2l}<\ff 1 2$ as $l<\ff d{(d-1)^+}$.
This together with \eqref{9} implies that
\beq\label{10} \|\mu_t\|_{1*}\le c(l,t)c_4(t) t^{-\ff{d(p-1)}{2p}}<\infty,\ \ \ t\in (0, \tau(\gg)).\end{equation}

Now, for any $T\in (s,\infty)$, let $(\mu_t)_{t\in [s,T]}$ be the solution to  \eqref{NFK} with initial value  $\mu_s=P_s^{\ast}\gamma$ at   time $s$ such that \eqref{*E} holds. When $b^{(1)}=0$ or  $\gg(|\cdot|)<\infty$, the estimate \eqref{MM} in Proposition \ref{L*1} implies
$$\E\int_s^T |b_t^{(1)}(X_t, \mu_t)|\d t<\infty.$$
Combining this   with $\|\si\|_\infty<\infty$, $b^{(i)}=0 \mbox{ with }2\leq i\leq l_0$ and \eqref{*E} as assumed,
we obtain
$$\int_s^T \mu_t\big(|b_t(\cdot,\mu_t)|+ \|\si_t\|^2\big)\d t<\infty.$$
Hence,  the superposition principle (see \cite{BR,T}) implies that the SDE
\beq\label{NNR}\d X_{s,t}^{\mu}= b_t(X_{s,t}^\mu, \mu_t)\d t+\si_t(X_{s,t}^\mu)\d W_t,\ \ t\in [s,T],\ \L_{X_{s,s}^\mu}=\mu_s\end{equation}
has a weak solution with $\L_{X_{s,t}^\mu}= \mu_t, t\in [s,T].$ Moreover,  by \eqref{*E}, Lemma \ref{L3'} holds for $l_0=2$ and $b_t^{(0)}(\cdot,\mu_t)$ in place of $b_t^{(2)}(\cdot,\mu_t)$, so that we derive \eqref{ES1'} and \eqref{ES2'} for $P_{t}^\mu$ in place of $\bar P_{s,t}^\mu$ as $b^{(i)}=0, 2\leq i\leq l_0$. Hence,
$$\sup_{t\in [s,T] } \|\mu_t\|_{k*} =\sup_{t\in [s,T] }\sup_{\|f\|_{\tt L^k}\le 1} |\mu_s(P_{s,t}^\mu f)|\le \|\mu_s\|_{k*} \sup_{t\in [s,T] }\|P_{s,t}^\mu\|_{\tt L^k\to\tt L^k} <\infty.$$
Combining this   weak solution of \eqref{NNR} with   the unique weak $\C_{p,k}$-solution of \eqref{E0} up to time $s$,  we may construct a   weak $\C_{p,k}$-solution for \eqref{E0} up to time $T$.
Therefore, by the above step (a),  \eqref{E0} has a unique weak $\C_{p,k}$-solution up to time $T$,  so that $\tau(\gg)> T$. Since $T\in (s,\infty)$ is arbitrary, we obtain $\tau(\gg)=\infty.$

Finally, by Theorem \ref{T0}(4), when $(A_2)$ holds, \eqref{E0}   has a unique global strong $\C_{p,k}$-solution.
Finally, repeating the proof of \eqref{9} for $q\in [1,p]$   replacing $1$, we prove \eqref{ER}.

 \end{proof}

\section{Proof of Theorem \ref{T02}}

 By Theorem \ref{T0},
for any  $\gamma\in\scr P_{p*}$ and $T\in (0,\tau(\gg))$,   we have
$$\mu :=(P_t^*\gg)_{t\in [0,T]}\in \C_{p,k}^{\gg,T}.$$
For simplicity, in the following we denote by $P_{s,t}^\gg$ the  operator  $P_{s,t}^\mu$ defined in \eqref{PM'} for   $\mu_t=P_t^\ast\gamma$, i.e.
\begin{align}\label{kgy}  P_{s,t}^\gamma f(x)=\E [f({X}_{s,t}^{\gg,x})],\ \ 0\le s\le t< \tau(\gamma),\ f\in\scr B_b(\R^d),\ x\in\R^d,
 \end{align}
where for fixed   $(s,x)\in [0,\tau(\gg))\times\R^d$,    $({X}_{s,t}^{\gg,x})_{  t\in [s, \tau(\gamma))}$ is the unique (weak) solution to  the SDE
$$\d  {X}_{s,t}^{\gg,x}=b_t({X}_{s,t}^{\gg,x}, P_t^\ast\gamma)\d t+\sigma_t({X}_{s,t}^{\gg,x})\d W_t,\ \ {X}_{s,s}^{\gg,x}=x,\ t\in [s,\tau(\gg)).$$
 Moreover, simply denote $P_t^\gg=P_{0,t}^\gg$ for $t\in [0,\tau(\gamma)).$

We first establish the estimates in  Lemma \ref{L3'} for $P_{s,t}^\gg$ in place of $\bar{P}_{s,t}^\mu$, which is crucial in the proof of Theorem \ref{T02}(1).

\beg{lem}\label{L01}  Assume  $(A_1)$ for   $b$   in $\eqref{B}$, and let $\kk_t(\gg)$ be in $\eqref{1ga}$.
Then for any $1<p_1\le p_2\le \infty$,  the following assertions hold for some  increasing  function $\bb: (0,\infty)\to (0,\infty).$
\beg{enumerate} \item[$(1)$] For any $\gg\in \scr P_{p*}$,
\beq\label{ES3} \|P_{s,t}^\gg\|_{\tt L^{p_1}\to\tt L^{p_2}}\le \bb_t \e^{\bb_t t^{2\theta} \kk_t(\gg)^2  } (t-s)^{-\ff{d(p_2-p_1)}{2p_1p_2}},
\ \ 0\le s<t< \tau(\gg). \end{equation}
Consequently, for any $n\in\mathbb N$ there exists a constant $c(n)\in (0,\infty)$ such that
\beq\label{ES3'} \|P_{s,t}^\gg\|_{\tt L^{p_1}\to\tt L^{p_2}}\le c(n)   (t-s)^{-\ff{d(p_2-p_1)}{2p_1p_2}},
\ \ 0\le s<t\le \tau_n(\gg),\ \gg\in \scr P_{p*}. \end{equation}
\item[$(2)$] If   $(A_2)$ holds,  then for any $\gg\in \scr P_{p*},$
\beq\label{ES4}   \|\nn P_{s,t}^\gg\|_{\tt L^{p_1}\to\tt L^{p_2}} \le \bb_t  \e^{\bb_t t\kk_t(\gg)^{\theta^{-1}}}  (t-s)^{-\frac{1}{2}-\ff{d(p_2-p_1)}{2p_1p_2}},
\ \    \ 0\le s<t< \tau(\gg).\end{equation}
Consequently, for any $n\in\mathbb N$, there exists a constant $c(n)\in (0,\infty)$ such that
\beq\label{ES4'} \|\nn P_{s,t}^\gg\|_{\tt L^{p_1}\to\tt L^{p_2}}\le c(n)   (t-s)^{-\frac{1}{2}-\ff{d(p_2-p_1)}{2p_1p_2}},\ \ 0\le s<t\le \tau_n(\gg),\ \gg\in \scr P_{p*}.\end{equation}
  \end{enumerate} \end{lem}

\beg{proof}  Without loss of generality, we only prove  for $s=0$ and   $t\in (0,\tau(\gg)).$
Let $\mu_s:=P_s^*\gg, s\in [0,t]$, and let $(\bar P_{s,s'}^\mu)_{0\le s\le s'\le t}$ be defined as in \eqref{SM}.
When $b^{(0)}=0$ we have $P_t^\gg=\bar P_t^\mu$, so that the desired estimates follow from \eqref{ES1'} and \eqref{ES2'}. It suffices to consider the case that $b^{(0)}\ne 0$.
Simply denote $\tt p_1=\ss{p_1}.$ In the following, all positive constants $\{c_i(t)\}_{i\ge 0}$ are  increasing in $t\in (0,\infty)$.

 (1) By $(A_1)$ and \eqref{EST}, we find    $c_0(t),c_1(t)\in (0,\infty)$ such that
 \beq\label{PO}|b_t^{(0)}(\cdot,\mu_t)|\leq c_0(t)\|\mu_t\|_{k*}\le c_1(t)\kk_t(\gg)t^{-\ff{d(p-k)}{2pk}}, \ \ t\in (0,\tau(\gg)).\end{equation} Since $(p,k)\in \D$ implies  $\theta:=\ff 1 2-\ff{d(p-k)}{2pk}>0$,
  by Girsanov's theorem and H\"older's inequality, we find   $c_2(t)\in (0,\infty)$  such that
\beq\label{EN1}|P_t^\gg f(x)|=|\E[R_tf(\bar X_{t}^{\mu,x})]|\le c_2(t)\e^{c_2(t)\kk_t(\gg)^2 t^{2\theta}}  (\bar P_{t}^\mu |f|^{\tt p_1})^{\ff 1 {\tt p_1}},\ \ t\in (0, \tau(\gg)),\end{equation}
where
$$R_t:=\e^{\int_0^t \<\eta_r,\d W_r\>-\ff 1 2\int_0^t |\eta_r|^2\d r},\ \ \eta_r:=(\si_r^{\ast}a_r^{-1}b_r^{(0)}(\cdot,\mu_r))(\bar X_{r}^{\mu,x}).$$
Combining \eqref{EN1} with \eqref{ES1'} and \eqref{1ga}, we find    $c_3(t)\in (0,\infty)$ such that
 \beg{align*}&  \|P_t^\gg \|_{\tt L^{p_1}\to\tt L^{p_2}}:= \sup_{\|f\|_{\tt L^{p_1}}\le 1} \|P_t^\gg f\|_{\tt L^{p_2}} \le \sup_{\|f\|_{\tt L^{p_1}}\le 1}
c_2(t)\e^{c_2(t)t^{2\theta}\kk_t(\gg)^2 }  \|\bar P_{t}^\mu|f|^{\tt p_1}\|_{\tt L^{p_2/\tt p_1}}^{1/\tt p_1}\\
&=c_2(t)\e^{c_2(t)t^{2\theta}\kk_t(\gg)^2  }  \|\bar P_{t}^\mu \|_{\tt L^{\tt p_1}\to\tt L^{ p_2/\tt p_1}}^{1/\tt p_1}  \le c_3(t)  \e^{c_2(t)t^{2\theta}\kk_t(\gg)^2 }  t^{-\ff{d(p_2-p_1)}{2p_1p_2}},\ \ \ t\in (0,\tau(\gg)).\end{align*}
So, \eqref{ES3} holds for some increasing function $\bb: (0,\infty)\to (0,\infty)$. Noting that \eqref{TT0} and \eqref{EST} imply
\beq\label{TT0*}t^{2\theta}\kk_t(\gg)^2\le ( \bb_1(n)\|\gg\|_{p*})^2  \bb_0(n)^{2\theta}\|\gg\|_{p*}^{-2}=  \bb_1(n)^2\bb_0(n)^{2\theta},\ \ t\le \tau_n(\gg),\end{equation}     \eqref{ES3'} follows from \eqref{ES3}.

(2)   By the same reason leading to \eqref{nab}, the estimates \eqref{DR} and the Bismut formula \eqref{BSM} in Proposition \ref{L*1} enable us to  find   $k_1(t,\gg)\in (0,\infty)$   such that
$$|\nn P_t^\gg f|\le k_1(t,\gg) t^{-\ff 1 2} (P_t^\gg |f|^{\tt p_1})^{1/\tt p_1},\ \ t\in (0,\tau(\gg)),\ \gg\in \scr P_{p*}.$$
Combining this  with \eqref{ES3} and the   argument deducing \eqref{ES2'} from \eqref{nab}, we find   $k_2(t,\gg)\in (0,\infty)$   such that
\beq\label{ES4''} \|\nn P_t^\gg\|_{\tt L^{p_1}\to\tt L^{p_2}}\le k_2(t,\gg)   t^{-\ff 1 2-\ff{d(p_2-p_1)}{2p_1p_2}},\  \ t\in (0,\tau(\gg)),\ \gg\in \scr P_{p*}.\end{equation}
To derive \eqref{ES4} with   $\bb_t\in (0,\infty)$ independent of $\gamma$, we apply   the Duhamel formula \eqref{DH1} for $P_t^\mu=P_t^\gg$ as $\mu_t=P_t^*\gg$,
which together with \eqref{PO},  \eqref{ES1'} and \eqref{ES2'}   implies that for some $c_4(t)\in (0,\infty)$
 \beq\label{DH3}\beg{split}  & \|\nabla P_t^\gamma\|_{\tilde{L}^{p_1}\to \tilde{L}^{p_1}} \leq  c_4(t)t^{-\frac{1}{2}}
 +c_4(t)\kk_t(\gg) \int_0^t\|\nabla P_s^\gamma\|_{\tilde{L}^{p_2}\to \tilde{L}^{p_2}} s^{-\ff{d(p-k)}{2pk}} (t-s)^{-\frac{1}{2}}\d s,\\
 &\qquad \gg\in \scr P_{p*},\ \ t\in (0,\tau(\gg)).\end{split}\end{equation}
  By \eqref{ES4''}  for $p_1=p_2$, for any $\lambda\geq 0$,  we have
\beq \label{HD1} H_t:=\sup_{s\in(0,t]}\e^{-\lambda s}s^{\frac{1}{2}}\|\nabla P_s^\gamma\|_{\tilde{L}^{p_2}\to \tilde{L}^{p_2}}<\infty.
\end{equation}
It follows from \eqref{DH3} that
\beq\label{DH3'}  H_t\leq  c_4(t)+c_4(t)\kk_t(\gg)H_t\sup_{s\in[0,t]}s^{\frac{1}{2}}\int_0^s r^{-\frac{1}{2}-\ff{d(p-k)}{2pk}} \e^{-\lambda(s-r)} (s-r)^{-\frac{1}{2}}\d r.
\end{equation}
By the FKG and H\"older inequalities, we can find a constant $c_5\in (0,\infty)$ such that
\beq\label{FKG}\begin{split}  &s^{\frac{1}{2}}\int_0^s r^{-\frac{1}{2}-\ff{d(p-k)}{2pk}}\e^{-\ll(s-r)} (s-r)^{-\ff 1 2}  \d r\\
&\le  s^{-\frac{1}{2}}\bigg(\int_0^s r^{-\frac{1}{2}-\ff{d(p-k)}{2pk}}\d r\bigg)\int_0^s\e^{-\ll(s-r)}  (s-r)^{-\ff 1 2}  \d r\\
 &=  \theta^{-1}s^{-\ff{d(p-k)}{2pk}}\bigg(\int_0^s\e^{-\ff{\ll}{\theta}(s-r)} \d r\bigg)^\theta \bigg(\int_0^s (s-r)^{-\ff 1 {2(1-\theta)}}\d r\bigg)^{1-\theta} \\
&\le c_5  \ll^{-\theta},\ \ \ll,s\in (0,\infty).\end{split}\end{equation}
Substituting this into \eqref{DH3'}, we conclude
$$H_t\leq  c_4(t)+c_4(t)c_5\kk_t(\gg)\lambda^{-\theta}H_t.$$
By  $H_t<\infty$ and taking $\lambda =[2c_4(t)c_5\kk_t(\gg)]^{\theta^{-1}}$, we get $H_t\leq 2c_4(t)$, which together with \eqref{HD1} yields that for some
$c_6(t)\in (0,\infty)$
\begin{align}\label{p-p}
\|\nabla  P^\gamma_{t}\|_{\tilde{L}^{p_2}\to \tilde{L}^{p_2}}\leq c_6(t)\e^{c_6(t)t \kk_t(\gg)^{\theta^{-1}} } t^{-\frac{1}{2}},\ \ t\in (0,\tau(\gg),\ \gg\in \scr P_{p*}.
\end{align}
  By \eqref{DH1} for $P_t^\mu=P_t^\gg$ since $\mu_t=P_t^\ast\gg$,   \eqref{p-p}, $(A_1)$, \eqref{PO} and \eqref{ES2'},     we find constants $K(t), c_7(t),c_8(t),c_9(t)\in (0,\infty)$ such that
  for any $\gg\in \scr P_{p*}$,
\begin{align*}
&\|\nabla P_t^\gamma\|_{\tilde{L}^{p_1}\to \tilde{L}^{p_2}} \\
&\le \|\nn \bar P_t^\mu\|_{\tilde{L}^{p_1}\to \tilde{L}^{p_2}}+K(t)\int_0^t \|\nn P_s^\gg\|_{\tt L^{p_2}\to\tt L^{p_2}}
\|P_s^*\gg\|_{k*} \|\nn\bar P_{s,t}^\mu \|_{\tt L^{p_1}\to\tt L^{p_2}}\d s\\
&\leq  c_7(t)t^{-\frac{1}{2}-\ff{d(p_2-p_1)}{2p_1p_2}}+c_7(t)\kk_t(\gg)\e^{c_6(t)t\kk_t(\gg)^{\theta^{-1}}}\int_0^t  s^{-\frac{1}{2}-\ff{d(p-k)}{2pk}} (t-s)^{-\frac{1}{2}-\ff{d(p_2-p_1)}{2p_1p_2}}\d s\\
&\leq c_8(t)t^{-\frac{1}{2}-\ff{d(p_2-p_1)}{2p_1p_2}}+c_8(t)\kk_t(\gg)\e^{c_6(t)t\kk_t(\gg)^{\theta^{-1}}} t^{-\ff{d(p-k)}{2pk}-\ff{d(p_2-p_1)}{2p_1p_2}}\\
&=c_8(t)\Big(1+t^{\theta} \kk_t(\gg) \e^{c_6(t)t\kk_t(\gg)^{\theta^{-1}}}\Big)t^{-\frac{1}{2}-\ff{d(p_2-p_1)}{2p_1p_2}}\\
&\le c_9 (t) \e^{c_9(t)t\kk_t(\gg)^{\theta^{-1}}} t^{-\frac{1}{2}-\ff{d(p_2-p_1)}{2p_1p_2}},\ \  \
t\in (0,\tau(\gg)).\end{align*}
Thus, \eqref{ES4} holds for some increasing $\bb: (0,\infty)\to (0,\infty)$, and it implies \eqref{ES4'} due to    \eqref{TT0*}.
 \end{proof}

 Combining Lemma \ref{L01} with  Proposition \ref{L*1} and Proposition \ref{L*2} addressed in Section 5,  we are ready to prove Theorem \ref{T02}.

\beg{proof}[\textbf{Proof of Theorem \ref{T02}(1)}] All constants $K$ and $c_i$ below may increasingly depend on $t>0$.   For fixed $\gg,\tt\gg\in \scr P_{p*}$, let $\pi\in \C(\gg,\tt\gg)$ such that
\beq\label{CPL} \W_q(\gg,\tt\gg)= \bigg(\int_{\R^d\times\R^d}|x-y|^q\pi(\d x,\d y)\bigg)^{\ff 1 q}.\end{equation}
For  $(P_{s,t}^\gamma)_{0\le s\le t\le T}$ defined in \eqref{kgy}, denote $P_t^\gg=  P_{0,t}^\gg$ and define
$P_t^{\gg*}:\scr P\to\scr P$ by
\beq\label{DF9} (P_t^{\gamma\ast}\nu)(f):=\int_{\R^d}P_t^\gamma f(x)\nu(\d x),\ \ f\in\scr B_b(\R^d), \ t\in [0,\tau(\gg)),\ \nu\in \scr P.\end{equation}
Let  $\hat{p}:=\ff{qp}{q-1}.$ By $(\hat p,k)\in\D$ implies
  $\ff{d(\hat{p}-k)}{\hat{p}k}\in [0,1)$.  By \eqref{ES4} we find a constant $c_1\in (0,\infty)$ such that
\beq\label{X1}   \|\nn P_t^\gg \|_{\tt L^k\to \tt L^{\hat p}} \le c_1  \e^{c_1t\kk_t(\gg)^{\theta^{-1}}}  t^{-\frac{1}{2}-\ff{d(\hat p-k)}{2\hat p k}},
   \ \ t\in (0,\tau(\gg)). \end{equation}
   Consider the maximal functional
\beq\label{MF}\scr M f(x):= \sup_{r\in (0,1)}\ff 1 {|B(x,r)|} \int_{B(x,r)} f(y)\d y,\ \ x\in\R^d,\end{equation}
   for a nonnegative measurable function $f$.
By \cite[Lemma 2.1]{XXZZ}  and  $P_t^\gg f\in C(\R^d),$   we find a constant $c_2\in (0,\infty)$ such that
  \beg{align*} &|P_t^\gg f(x)-P_t^\gg f(y)|\le c_2 |x-y|\big(\scr M|\nn P_t^\gg f|(x)+\scr M|\nn P_t^\gg f|(y)+\|P_t^\gg f\|_\infty\big),\\
 & \big\|\scr M|\nn P_t^\gg f|\big\|_{\tt L^{\ff{pq}{q-1}}}\le c_2 \|\nn P_t^\gg f\|_{\tt L^{\ff{pq}{q-1}}}, \ \ t\in (0,\tau(\gg)),\ x,y\in\R^d.
\end{align*}
 Combining this with H\"older's inequality, we find a constant $c_3\in (1,\infty)$ such that
\beg{align*} & \big\| P_t^*\gg- P_t^{\gg *}\tt\gg\big\|_{k*}=  \big\|P_t^{\gg *}\gg- P_t^{\gg *}\tt\gg\big\|_{k*}=\sup_{\|f\|_{\tt L^k}\le 1} \big|\gg(P_t^\gg f)-\tt\gg(P_t^\gg f)\big|\\
&=\sup_{\|f\|_{\tt L^k}\le 1} \bigg|\int_{\R^d\times\R^d} \big(P_t^\gg f(x)- P_t^\gg f(y)\big)\pi(\d x,\d y)\bigg|\\
&\le c_2\sup_{\|f\|_{\tt L^k}\le 1} \bigg|\int_{\R^d\times\R^d} |x-y|\big(\scr M|\nn P_t^\gg f|(x)+\scr M|\nn P_t^\gg f|(y)+\|P_t^\gg f\|_\infty\big)\pi(\d x,\d y)\bigg|\\
&\le  c_3  \W_q(\gg,\tt\gg) \sup_{\|f\|_{\tt L^k}\le 1}\Big[(\gg+\tt\gg)\big((\scr M|\nn P_t^\gg f|)^{\ff{q}{q-1}}\big)\Big]^{\ff {q-1}q}+c_2\W_q(\gg,\tt\gg)\|P_t^\gamma\|_{\tt L^k\to\tt L^{\infty}}\\
&\le  c_3  \W_q(\gg,\tt\gg) \sup_{\|f\|_{\tt L^k}\le 1} (\|\gg\|_{p*}+\|\tt\gg\|_{p*})^{\ff{q-1}q} \big\|\scr M|\nn P_t^\gg f|\big\|_{\tt L^{\ff{pq}{q-1}}}+c_2\W_q(\gg,\tt\gg)\|P_t^\gamma\|_{\tt L^k\to\tt L^{\infty}}\\
&\le c_2c_3 \W_q(\gg,\tt\gg)   (\|\gg\|_{p*}+\|\tt\gg\|_{p*})^{\ff{q-1}q} \left( \big\| \nn P_t^\gg \big\|_{\tt L^k\to\tt L^{\hat p}}+\|P_t^\gamma\|_{\tt L^k\to\tt L^{\infty}}\right),\ \  t\in (0,\tau(\gg)\land \tau(\tt\gg)).\end{align*}
This together with \eqref{X1} and \eqref{ES3} yields
\beq\label{X3} \beg{split} \big\|P_t^*\gg- P_t^{\gg *}\tt\gg\big\|_{k*}
 \le &\,c_1 c_2c_3   (\|\gg\|_{p*}+\|\tt\gg\|_{p*})^{\ff{q-1}q}\e^{c_1t\kk_t(\gg)^{\theta^{-1}}}  t^{-\left(\left(\frac{1}{2}+\ff{d(\hat{p}-k)}{2\hat{p}k}\right)\vee(\frac{d}{2k})\right)}\W_q(\gg,\tt\gg),\\
 &   \gg,\tt\gg\in \scr P_{p*},\ t\in (0,\tau(\gg)\land \tau(\tt\gg)).\end{split}\end{equation}

On the other hand, by Duhamel's formula  \eqref{DH1'} below, we have
$$P_t^{\gg}f-P_t^{\tt\gg}f= \int_0^t P_s^\gg \big\<b_s(\cdot,P_s^*\gg)- b_s(\cdot,P_s^*\tt\gg),\ \nn P_{s,t}^{\tt\gg}f\big\>\d s,\ \ f\in C_0^\infty(\R^d),$$
and   $(A_1)$ implies
$$ |b_t(x,P_t^*\gg)-b_t(x,P_t^*\tt\gg)|\le K\|P_t^*\gg-P_t^*\tt\gg\|_{k*},\ \ t\in [0, \tau(\gg)\land  \tau(\tt\gg)].$$
 Then
\beg{align*}
&\big|\big(P_t^{\gg *} \tt\gg\big)(f)- \big(P_t^*\tt\gg\big)(f)\big|=\big|(P_t^{\gg *}\tt\gg)(f)
 - (P_r^{\tt\gg *}\tt\gg)(f)\big|
=\big|\tt\gg\big(P_t^{\gg}f-P_t^{\tt\gg}f\big)\big|\\
&\le \int_0^t \tt\gg\Big(P_s^\gg\big(\big|b_s(\cdot,P_s^*\gg)- b_s(\cdot,P_s^*\tt\gg)\big|\cdot\big|\nn P_{s,t}^{\tt\gg}f\big|\big)\Big)\d s\\
&\le K \|\tt\gg\|_{p*} \int_0^t \big\|P_s^*\gg-P_s^*\tt\gg\big\|_{k*}  \|P_s^\gg\|_{\tt L^p\to\tt L^p}
\|\nn P_{s,t}^{\tt\gg}f\|_{\tt L^p}\d s, \ \  t\in (0,\tau(\gg)\land \tau(\tt\gg)).\end{align*}
This together with  \eqref{ES3} and \eqref{ES4} yields that for some constant $c_4\in (0,\infty)$
\beg{align*}\big\| P_t^{\gg *}\tt\gg-P_t^*\tt\gg\big\|_{k*}
&\le c_4 \|\gamma\|_{p*}\e^{c_4t^{2\theta}\kk_t(\gg)^2+c_4t \kk_t(\tt\gg)^{\theta^{-1}}}
 \int_0^t \big\|P_s^*\gg-P_s^*\tt\gg\big\|_{k*} (t-s)^{-\ff 1 2-\ff{d(p-k)}{2pk}} \d s,\\
 &\qquad \ \gg,\tt\gg\in \scr P_{p*}, \  t\in [0, \tau(\gg)\land  \tau(\tt\gg)).\end{align*}
Combining this with \eqref{X3} and the triangle inequality, we obtain
\beg{align*}  &\big\|P_t^*\gg-P_t^*\tt\gg\big\|_{k*}\le \big\|P_t^*\gg- P_t^{\gg *}\tt\gg\big\|_{k*}
+\big\| P_t^{\gg *}\tt\gg-P_t^*\tt\gg\big\|_{k*}\\
&\le  c_1 c_2c_3    (\|\gg\|_{p*}+\|\tt\gg\|_{p*})^{\ff{q-1}q} \e^{c_1t\kk_t(\gg)^{\theta^{-1}}}t^{-\left(\left(\frac{1}{2}+\ff{d(\hat{p}-k)}{2\hat{p}k}\right)\vee(\frac{d}{2k})\right)}
\W_q(\gg,\tt\gg)    \\
&\quad +c_4  \|\gamma\|_{p*}\e^{c_4t^{2\theta}\kk_t(\gg)^2+c_4t\kk_t(\tt\gg)^{\theta^{-1}}}   \int_0^t \big\|P_s^*\gg-P_s^*\tt\gg\big\|_{k*} (t-s)^{-\ff 1 2-\ff{d(p-k)}{2pk}} \d s,\\
&\qquad \gg,\tt\gg\in \scr P_{p*}, \  t\in [0, \tau(\gg)\land  \tau(\tt\gg)).\end{align*}
 Note that $\ff{d(\hat{p}-k)}{2\hat{p}k}=\ff{d(pq-(q-1)k)}{2pqk}$.  So,  for any constant $\ll\in (0,\infty)$ and  $t\in [0, \tau(\gg)\land  \tau(\tt\gg)),$  the  finite quantity
\beq\label{YU}G_t:=\sup_{s\in (0, t]} \e^{-\ll s} s^{(\ff 1 2+\ff{d(pq-(q-1)k)}{2pqk})\vee\frac{d}{2k}}\big\|P_s^*\gg-P_s^*\tt\gg\big\|_{k*}  \end{equation}
 satisfies
 \begin{align*} G_t \le &\,c_1c_2c_3 (\|\gg\|_{p*}+\|\tt\gg\|_{p*})^{\ff{q-1}q}  \e^{c_1t\kk_t(\gg)^{\theta^{-1}}} \W_q(\gg,\tt\gg) \\
&+  c_4\|\gamma\|_{p*}  \e^{c_4t^{2\theta}\kk_t(\gg)^2+ c_4t\kk_t(\tt\gg)^{\theta^{-1}}} G_t \\
&\quad\times    \sup_{s\in (0, t]}s^{\ff 1 2+\ff{d(pq-(q-1)k)}{2pqk}}
\int_0^s r^{-\ff 1 2-\ff{d(pq-(q-1)k)}{2pqk}}\e^{-\lambda(s-r)}(s-r)^{-\ff 1 2-\ff{d(p-k)}{2pk}}\d r.
\end{align*}
Similarly to  \eqref{FKG},   we find a constant $c_5\in (0,\infty)$ such that
\beg{align*}&s^{\left(\left(\frac{1}{2}+\ff{d(\hat{p}-k)}{2\hat{p}k}\right)\vee(\frac{d}{2k})\right)}
\int_0^s r^{-\left(\left(\frac{1}{2}+\ff{d(\hat{p}-k)}{2\hat{p}k}\right)\vee(\frac{d}{2k})\right)}\e^{-\lambda(s-r)}(s-r)^{-\ff 1 2-\ff{d(p-k)}{2pk}}\d r\\
&\le s^{\left(\left(\frac{1}{2}+\ff{d(\hat{p}-k)}{2\hat{p}k}\right)\vee(\frac{d}{2k})\right)}\bigg(\ff 1 s\int_0^s r^{-\left(\left(\frac{1}{2}+\ff{d(\hat{p}-k)}{2\hat{p}k}\right)\vee(\frac{d}{2k})\right)}\d r\bigg)\int_0^s\e^{-\lambda(s-r)}(s-r)^{-\ff 1 2-\ff{d(p-k)}{2pk}} \d r\\
&\le c_5  \ll^{-\theta},\ \ \ \theta:= \ff 1 2-\ff{d(p-k)}{2pk}>0,\end{align*}
so that we obtain
\beg{align*}  \Big(1- c_4c_5  \|\gamma\|_{p*}\e^{c_4t^{2\theta}\kk_t(\gg)^2+c_4t\kk_t(\tt\gg)^{\theta^{-1}}} \ll^{-\theta} \Big)G_t \le  c_1c_2c_3 (\|\gg\|_{p*}+\|\tt\gg\|_{p*})^{\ff{q-1}q}   \e^{c_1t\kk_t(\gg)^{\theta^{-1}}} \W_q(\gg,\tt\gg).\end{align*}
Taking
$$\ll= \Big[2c_4c_5 \|\gamma\|_{p*} \e^{c_4t^{2\theta}\kk_t(\gg)^2+c_4t\kk_t(\tt\gg)^{\theta^{-1}}} \Big]^{\theta^{-1}},$$
we derive
$$G_t\le 2c_1c_2c_3 (\|\gg\|_{p*}+\|\tt\gg\|_{p*})^{\ff {q-1}q} \e^{c_1t\kk_t(\gg)^{\theta^{-1}}} \W_q(\gg,\tt\gg),$$
which together with the definition of $G_t$ in \eqref{YU},  $$t^{2\theta}\kk_t(\gg)^2=\left(t\kk_t(\gg)^{\theta^{-1}}\right)^{2\theta}\leq 1+t\kk_t(\gg)^{\theta^{-1}},$$ and $\kk_t(\gg)\geq \|\gamma\|_{p*}$ due to \eqref{1ga}  implies that for some constant $c_6\in (0,\infty),$
\begin{align*}\|P_t^*\gg-P_t^*\tt\gg\|_{k*}&\le  (\|\gg\|_{p*}+\|\tt\gg\|_{p*})^{\ff{q-1}q}  t^{-\left(\left(\frac{1}{2}+\ff{d(\hat{p}-k)}{2\hat{p}k}\right)\vee(\frac{d}{2k})\right)} \exp\left[c_6\e^{c_6t\kk_t(\gg)^{\theta^{-1}}+c_6t\kk_t(\tt\gg)^{\theta^{-1}}} \right].
\end{align*}
This implies \eqref{ES5} for some $\bb: (0,\infty)\to (0,\infty)$.

Finally, when $p=\infty$ or $b^{(0)}=0$, $\kk_t(\gg)$ defined in \eqref{1ga} is bounded above by some constant $c(t)\in (0,\infty)$ uniformly in $\gg\in \scr P_{p*}$. So,
\eqref{ES5} implies  \eqref{ES5'}.
\end{proof}

\beg{proof}[\textbf{Proof of Theorem \ref{T02}(2)}]  For fixed $t\in [0,\tau(\gg)\land \tau(\tt\gg)),$ denote
  $$\gg_t:=P_t^*\gg,\ \ \ \tt\gg_t:= P_t^*\tt\gg.  $$
  To estimate the relative entropy $\Ent(\gg_t|\tt\gg_t)$,   we consider the SDEs
  \beq\label{LG} \beg{split} &\d X_s= b_s(X_s,\gg_s)\d s +\si_s(X_s)\d W_s,\\
  &\d Y_s= b_s(Y_s,\tilde{\gg}_s)\d s +\si_s(Y_s)\d W_s,\ \ \ s\in [0,t],\end{split}\end{equation}
  such that the initial values $X_0,Y_0$ are $\F_0$-measurable   satisfying
  \beq\label{CP} \L_{X_0}=\gg,\ \ \  \L_{Y_0}=\tt\gg,\ \ \ \E|X_0-Y_0|^2= \W_2(\gg,\tt\gg)^2.\end{equation}
 Note that we can always choose suitable $\F_0$ independent of $W_t$ such that the above  $X_0$ and $Y_0$ exist.
  Then
\beq\label{GT}  \gg_t=P_t^*\gg=\L_{X_t},\ \ \ \tt\gg_t=P_t^*\tt\gg=\L_{Y_t}.\end{equation}
  By $(A_1)$ and  \eqref{ES5},  we find a constant $K(t)\in (0,\infty)$  increasing in $t$   such  that
\beq\label{XX0}\xi_s := (\si_s^*a_s^{-1})(Y_s)  \big[b_s(Y_s, \gg_s)- b_s(Y_s, \tt\gg_s)\big]\end{equation} satisfies
\beq\label{XX}|\xi_s |^2 \le K(t)(\|\gg\|_{p*}\lor\|\tt\gg\|_{p*})^{\ff{2(q-1)}q}K_{t,\bb}^{p,k}(\gg,\tt\gg)^2 \W_q(\gg,\tt\gg)^2 s^{-\left((1+\ff{d(qp-(q-1)k)}{pqk})\vee \frac{d}{k}\right)},\ \ s\in (0,t].\end{equation}
 Since
\beq\label{SG} \int_0^t s^{-\left((1+\ff{d(qp-(q-1)k)}{pqk})\vee \frac{d}{k}\right)}\d s=\infty,\end{equation}
we can not apply  Girsanov's theorem to kill $\xi_s$.
To overcome the singularity of $|\xi_s|^2$ for small $s>0$,     we will
  apply the bi-coupling argument developed in \cite{23RW}, and finish the proof in the following  three steps.

 (a)
 We first establish the   log-Harnack inequality  for $P_{t}^\gg:$ for any $\theta'\in (0,\theta)$,  there exists   $c_0 (t)  \in (0,\infty)$ increasingly in $t$  such that
\beq\label{cty}\beg{split}
&P_{s,t}^\gamma \log f(x) \leq \log P_{s,t}^\gamma f(y)+\frac{c_0(t) |x-y|^2}{s_t(\theta',\gg)\wedge (t-s)},\\
&\qquad  \ x,y\in\R^d,   0\leq s< t<\tau(\gg),\    \gg\in \scr P_{p*},\  f\in \B_b^+(\R^d)
\end{split}\end{equation}
 for $s_t(\theta',\gg)$ defined in \eqref{ST}.
 We will prove this estimate by applying Proposition \ref{L*1}(4) to
\beq\label{B00}b^{0,1}_t:= b_t^{(0)}(\cdot,\gg_t),\ \ \ b_t^{0,i}:= b^{(i)}_t(\cdot,\gg_t),\ 2\le i\le l_0.\end{equation}
 By $(A_1)$, we have
 \beq\label{B2} \sup_{\gg\in \scr P_{p*}}\sup_{2\le i\le l_0}\|b^{0,i}\|_{\tt L_{q_i'}^{p_i'}(s_t(\theta',\gg))}<\infty,\ \ t\in (0,\tau(\gamma)).\end{equation}
 Next,  $\theta'\in (0,\theta)$ implies
 $$q_1':= \Big(\ff{d(p-k)}{2pk}+\theta'\Big)^{-1}\in \Big(2,\ff{2pk}{d(p-k)}\Big).$$
 Then   there exists $p_1'\in (d,\infty)$ such that  $(p_1',q_1')\in \scr K$. By
 \eqref{PO} we find constants $k_1,k_2\in (0,\infty)$ such that $b_t^{0,1}:=b^{(0)}_t(\cdot,\gg_t)$ satisfies
\beq\label{GM} \|b^{0,1}\|_{\tt L_{q_1'}^{p_1'}(s,t)}\le k_1c_1(t) \kk_t(\gg)  \bigg(\int_s^t r^{-\ff{q_1'd(p-k)}{2pk}}\d r\bigg)^{\ff 1 {q_1'}}
 \le k_2 c_1(t) \kk_t(\gg)(t-s)^{\theta'},\end{equation} where $c_1(t)\in (0,\infty)$ is increasing in $t$.
 By \eqref{ST},  we find a constant $k_3\in (0,\infty)$ such that
 $$\kk_t(\gg)(t-s)^{\theta'}\le k_3,\ \ 0<t-s\le  s_t(\theta', \gg).$$
 This together with \eqref{GM} implies that for a constant $k_4 \in (0,\infty)$   such that
  $$\|b^{0,1}\|_{\tt L_{q_1'}^{p_1'}(s,t)}\le k_4c_1(t),\ \ 0<t-s\le s_t(\theta',\gg), \ t\in (0,\tau(\gg)).$$
Combining this with \eqref{B2}, we may apply Proposition \ref{L*1}(4) to find
     $k_5\in (0,\infty)$  such that for any $f\in \B_b^+(\R^d)$,
\beq\label{o1}  \beg{split}
 & P_{s,t}^\gg \log f(x)\le \log P_{s,t}^\gg f(y) + \ff{k_5c_1(t)|x-y|^2}{t-s}, \\
  &\qquad \ x,y\in\R^d,\ \ \  0<t-s\le s_t(\theta',\gg),\ t\in (0,  \tau(\gg)). \end{split}\end{equation}
 Hence, \eqref{cty} holds for $t-s \le s_t(\theta',\gg)$ and $c_0(t)= k_5c_1(t).$

 Now let   $t-s>s_t(\theta',\gg),\ t\in (0,\tau(\gg)). $  By the semigroup property and Jensen's inequality,  we deduce from \eqref{o1} that
  \beg{align*}&
 P_{s,t}^\gg \log f(x)=P_{s,s+s_t(\theta',\gg)}^\gg P_{s+s_t(\theta',\gg),t}^\gg \log f(x)\leq P_{s,s+s_t(\theta',\gg)}^\gg \log P_{s+s_t(\theta',\gg),t}^\gg  f(x)\\
 &\leq \log P_{s,s+s_t(\theta',\gg)}^\gg P_{s+s_t(\theta',\gg),t}^\gg  f(y)+\ff{k_5c_1(t)|x-y|^2}{s_t(\theta',\gg)}\\
  &= \log P_{s,t}^\gg f(y)+ \ff{k_5c_1(t)|x-y|^2}{s_t(\theta',\gg)}, \ \ x,y\in\R^d,\ f\in \B_b^+(\R^d).\end{align*}
 So, \eqref{cty} also  holds for  $t-s>s_t(\theta',\gg)$ and $c_0(t):= k_5c_1(t).$

 (b) To apply the bi-coupling argument, for
 fixed $t\in (0,\tau(\gamma)\wedge \tau(\tilde{\gamma})),$  let
\begin{align}\label{key}
t'=\frac{t}{2}\wedge  s_t(\theta',\gamma)\wedge s_t(\theta',\tilde{\gamma}).
\end{align}
To cancel the singularity in \eqref{SG} for small $s>0$,  we construct the following SDE which will be coupled with two SDEs in \eqref{LG} respectively:
 \beq\label{BQ}   \d Z_s =  \big\{1_{[0,t']}(s) b_s(Z_s, \tt \gg_s)+ 1_{(t',t]}(s) b_s(Z_s,\gg_s)\big\}\d s+ \si_s(Z_s)\d W_s,\ \ Z_0=Y_0,s\in [0,t].\end{equation}
 By \eqref{GT} and     \cite[Lemma 2.1]{23RW},   we have
\beq\label{RW} \beg{split}
&\Ent(\gg_t| \tt\gg_t)=\Ent(\L_{X_t}|\L_{Y_t})\\
&\le 2 \Ent(\L_{X_t}|\L_{Z_t}) +  \log \int_{\R^d} \Big(\ff{\d\L_{Z_t}}{\d\L_{Y_t}}\Big)^{2}\d \L_{Y_t}
  =:2 I_1+  I_2.\end{split}\end{equation}
 Below we estimate $I_1$ and $I_2$  respectively.

(i) Estimate $I_1$.
 Let  $X_{t', s}^x$ solve the SDE
$$\d X_{t',s}^x= b_s(X_{t',s}^x, \gg_s)\d s +\si_s(X_{t',s}^x)\d W_s,\ \ X_{t',t'}^x=x,\ s\in [t',t],$$
and define
$$P_{t',t}^\gg f(x):=\E[f(X_{t',t}^x)],\ \ \ f\in \B_b(\R^d),\ x\in\R^d.$$
By the Markov property, we have
\beq\label{OG}\E[f(X_t)]= \E[ (P_{t',t}^\gg f)(X_{t'})],\ \ \
\E[f(Z_t)] = \E[ (P_{t',t}^\gg f)(Z_{t'})].\end{equation}
This together with \eqref{cty} for $s=t'$ and Jensen's inequality implies
\beq\label{OG3}\E[\log f(X_t)]\le \log \E[f(Z_t)]+ \ff {2c_0(t)}{  s_t(\theta',\gg)}\E[|X_{t'}- Z_{t'}|^2],\ f\in \B_b^+(\R^d).\end{equation}
By $(A_1)$ and $\gg,\tt\gg\in \C_{p,k}^t$,
$\tt b^{0,1}_s:=b_s(\cdot,\tt\gamma_s)- b_s(\cdot,\gg_s)$ satisfies
$\|\tt b^{0,1}\|_{\tt L_{q_1'}^{p_1'}(t)}<\infty.$ For $b^{0,1}$ and $b^{0,2}$ in \eqref{B00}, we have
\beg{align*} &b_s(\cdot,\gg_s) =  \sum_{i=1}^{l_0}b_s^{0,i}+ b_s^{(1)}(\cdot,\gg_s),\\
 &   b_s(\cdot,\tt\gg_s)  = \tt b_s^{0,1} +\sum_{i=1}^{l_0}b_s^{0,i}   + b_s^{(1)}(\cdot,\gg_s),\ \ \ s\in [0,t],\ x\in\R^d.
\end{align*} By $(A_1)$,  \eqref{PO}, \eqref{ST}  and $t'\le s_t(\theta',\gg)\land s_t(\theta',\tt\gg),$
we find $c_2(t)\in (0,\infty)$ increasing in $t\in (0,\infty)$ such that
\beg{align*}&\sum_{i=1}^{l_0}\|b^{0,i}\|_{\tt L_{q_i'}^{p_i'}(t')}+ \|\tt b^{0,1}\|_{\tt L_{q_1'}^{p_1'}(t')} \le c_2(t),\\
&\|b_s(\cdot,\gg_s)-b_s(\cdot,\tt\gg_s)\|_\infty\le c_2(t) \|\gg_s-\tt\gg_s\|_{k*},\ \ s\in [0,t'].\end{align*}  So,
by Proposition \ref{L*2}, we find $c_3(t)\in (0,\infty)$ increasing in $t\in (0,\infty)$ such that
\beq\label{MS} \E[|X_{t'}-Z_{t'}|^2]\le c_3(t)  \E|X_0-Y_0|^2+  c_3(t) \bigg(  \int_0^{t'} \|\gg_s-\tt\gg_s\|_{k*}\d s \bigg)^2.\end{equation}
Moreover, by
 $$t'\le s_t(\theta',\gg)\land s_t(\theta',\tt\gg),\ \  t'\kk_{t'}(\gg)^{1/\theta}=(t')^{1-\frac{\theta'}{\theta}}(t')^{\frac{\theta'}{\theta}}\kk_{t'}(\gg)^{1/\theta},$$
  \eqref{Ka} and \eqref{ST}, we find a constant $c_4(t)\in (0,\infty)$ increasing in $t\in (0,\infty) $ such that
 $$K_{t',\bb}^{p,k}(\gg,\tt\gg)\le c_4(t).$$
   Combining this with \eqref{ES5} and letting 
   $\tilde{\theta}=1-\left((\frac{1}{2}+\ff{d(qp-(q-1)k)}{2pqk})\vee \frac{d}{2k}\right)$, we find
 $c_5(t)\in (0,\infty)$ increasing in $t\in (0,\infty)$ such that
\beq\label{Xi}\beg{split}&\int_0^{t'} \|\gg_s - \tt\gg_s\|_{k*}\d s\le c_5(t) (t')^{\tilde{\theta}}  (\|\gg\|_{p*}\lor\|\tt\gg\|_{p*})^{\ff{q-1}q}   \W_q(\gg,\tt\gg) ,\ \  t\in[0,\tau(\gamma)\wedge \tau(\tilde{\gamma})).  \end{split}\end{equation} Combining this with \eqref{CP}, \eqref{MS} and $\|\gamma\|_{p*}\geq \|\gamma\|_{\infty*}=1$,  we find   $c_6(t)\in (0,\infty)$ increasing in $t\in (0,\infty)$ such that
\beg{align*} \E[|X_{t'}- Z_{t'}|^2] &\le c_6(t) (\|\gg\|_{p*}\lor\|\tt\gg\|_{p*})^{\ff{2(q-1)}q}  \Big( (t')^{2{\tilde{\theta}}}\W_q(\gg,\tt\gg)^2+\W_2(\gg,\tt\gg)^2\Big).\end{align*}
 Combining this with \eqref{key}, \eqref{OG3},     and the formula
$$\Ent(\mu|\nu)=\sup_{f\in \B_b^+(\R^d)} \big\{\mu(\log f)-\log \nu(f)\big\},\ \ \mu,\nu\in \scr P,$$
we find a constant $c_7(t) \in (0,\infty)$ increasing in $t\in (0,\infty)$ such that
 \beq\label{I11}\begin{split}
 &I_1:=\Ent(\L_{X_t}|\L_{Z_t})=\sup_{f\in \B_b^+(\R^d)} \big\{\E[\log   f(X_t)]-\log \E[ f(Z_t)]\big\},  \\
 & \leq   c_7(t)   (\|\gg\|_{p*}\lor\|\tt\gg\|_{p*})^{\ff{2(q-1)}q}\bigg( \ff{\W_2(\gg,\tt\gg)^2}{ s_t(\theta',\gamma)}+\ff{\W_q(\gg,\tt\gg)^2}{s_t(\theta',\gg)^{1-2\tilde{\theta}}}\bigg),
 \end{split}\end{equation}
where in the last step we have used $t'\le   s_t(\theta',\gg)$.

(ii) Estimate $I_2$. By \eqref{XX} for $\xi_s$   in \eqref{XX0},
 $$ R_s:=\e^{\int_{t'}^s \<\xi_r,\d W_r\>-\ff 1 2 \int_{t'}^s|\xi_r|^2\d r},\ \ s\in [t',t]$$
is a martingale, and by Girsanov's theorem,
$$\ff{\d\L_{Z_t}}{\d\L_{Y_t}}(Y_{t})= \E(R_t|Y_{t}).$$
Combining this with Jensen's inequality and \eqref{XX}, and denoting
$$\theta_t:= (\|\gg\|_{p*}\lor\|\tt\gg\|_{p*})^{\ff{2(q-1)}q}  K_{t,\bb}^{p,k}(\gg,\tt\gg)^2, $$
we find constants $c_9, c_{10}\in (0,\infty)$   such that
   \beg{align*}I_2&:= \log \E\bigg[\Big(\ff{\d\L_{Z_t}}{\d\L_{Y_t}}(Y_{t})\Big)^{2}\bigg]
\leq \log\E\Big[R_t^{2}\Big]\\
&\le \log\E\bigg[\e^{2\int_{t'}^t \<\xi_s,\d W_s\>- 2 \int_{t'}^t |\xi_s|^2\d s+   \big(\theta_t\int_{t'}^t s^{-\left(\left(1+\ff{d(qp-(q-1)k)}{pqk}\right)\vee\frac{d}{k}\right)}\d s\big)\W_q(\gg,\tt\gg)^2}\bigg]\\
&=  \bigg(\theta_t\int_{t'}^t s^{-\left(\left(1+\ff{d(qp-(q-1)k)}{pqk}\right)\vee\frac{d}{k}\right)}\d s\bigg)\W_q(\gg,\tt\gg)^2\\
&\le c_9   \theta_t \left((t')^{1-\left(\left(1+\ff{d(qp-(q-1)k)}{pqk}\right)\vee\frac{d}{k}\right)}-t^{1-\left( \left(1+\ff{d(qp-(q-1)k)}{pqk}\right)\vee\frac{d}{k}\right)}\right)\W_q(\gg,\tt\gg)^2\\
&\leq c_{10} \theta_t \left[s_t(\theta',\gamma)\land  s_t(\theta',\tilde\gamma)\right]^{1-\left(\left(1+\ff{d(qp-(q-1)k)}{pqk}\right)\vee\frac{d}{k}\right)}\W_q(\gg,\tt\gg)^2.
\end{align*}
By combining this with   \eqref{RW} and \eqref{I11}, we obtain \eqref{ES6} for some    $\bb: (0,\infty)\to (0,\infty).$

(c) If $p=\infty$, we have   $\scr P_{p*}=\scr P$, $\tau(\gg)=\infty$ and $\|\gg\|_{p*}= 1$ for any $\gg\in \scr P$, and we may take $q=1$ so that  $(\frac{pq}{q-1},k)=(\infty,k)\in \D.$   Hence,   \eqref{ES6} implies \eqref{ES7}.

If $b^{(0)}=0$, then $s_t(\theta',\gg)=t$ and $\kk_t(\gg)=0$  for $\gg\in \scr P_{p*}$, so that
  \eqref{ES8}  follows from \eqref{ES6}, \eqref{1ga} and \eqref{Ka}, for some different increasing function $\beta: [0,\infty)\to (0,\infty).$   \end{proof}

\section{  SDEs with several singular drifts}

In this part we present  some results on SDEs with several singular drifts,  which    include
well-posendness, regularity, the local hyperbound estimates on diffusion semigroup,  and Duhamel's formula. These results are used in the proofs of Theorem \ref{T0} and Theorem \ref{T02}, and extend
the existing ones   for SDEs with unique  singular drift.

  \subsection{The model and well-posedness}

We consider measurable   drifts
$$  b^{0,i},   b^{(1)}: [0,T]\times\R^d\to \R^d,\ \ 1\le i\le \ell',$$ where $T\in (0,\infty)$ and $\ell'\in\mathbb N$ are fixed. These drifts and
 $a:=\si\si^*$ satisfy the following assumption.
\beg{enumerate} \item[$(C)$] Let   $a:=\si\si^*$ satisfy
\beq\label{conti}\zeta(\vv):= \sup_{|x-y|\le\vv, t\in [0,T]} \|a_t(x)-a_t(y)\|\downarrow 0\ \text{as}\ \vv\downarrow 0.\end{equation}
There exist $ K\in (0,\infty)$ and $\{(p_i',q_i')\}_{1\le i\le \ell'} \subset \scr K $ such that
\beg{align*}
&\|  b^{0,i}\|_{\tt L_{q_i'}^{p_i'}(T)}+  \|a\|_\infty+\|a^{-1}\|_\infty+\|  b^{(1)}(0)\|_\infty\le K,  \\
& |  b^{(1)}_t(x)-  b_t^{(1)}(y)|\le K|x-y|,\ \ x,y\in \R^d,\ t\in [0,T].\end{align*}
\end{enumerate}
  For fixed $(s,x)\in [0,T)\times\R^d$, we consider the SDE
\beq\label{S0}\d   X_{s,t}^x= \bigg(\sum_{i=1}^{\ell'}  b_t^{0,i}+   b_t^{(1)}\bigg)(X_{s,t}^x)\d t +\si_t(X_{s,t}^x)\d W_t,\ \ t\in [s,T],\ X_{s,s}^x=x.\end{equation}
Simply denote $X_{t}^x=X_{0,t}^x.$  When the SDE \eqref{S0} is weakly well-posed, we define
$$P_{s,t} f(x):=\E\big[f(X_{s,t}^x)\big],\ \ 0\le s\le t\le T,\ x\in\R^d,\  f\in \B_b(\R^d).$$

\beg{prp}\label{P01} Assume $(C)$.
 Then for any  $(s,x)\in [0,T)\times\R^d$,    $\eqref{S0}$ is weakly well-posed.
If $(A_2)$ holds, then the SDE is strongly well-posed.
\end{prp}

\beg{proof}   According to \cite[Theorem 1.3.1]{RW24}, the assertions hold for $\ell'=1.$
Assume that the assertions hold for $\ell'=n$ for some $n\in\mathbb N$, it suffices to
prove for $\ell'=n+1$.

Let
\beg{align*} L_t = \ff 1 2{\rm tr}\{ a_t\nn^2\}+ \big\{  b^{0,1}_t+  b^{(1)}_t\big\}\cdot\nn, \
\ \ t\in [0,T].\end{align*}
By \cite[Theorem 2.1]{YZ},
 for large enough  $\ll\in (0,\infty)$, the PDE
\beq\label{PL0} (\pp_t +L_t-\ll)u_t=-  b_t^{0,1},\ \ \ t\in [0,T],\ u_T=0\end{equation}
has a unique solution $u: [0,T]\times \R^d\to\R^d$ such that for any  $\theta\in (1,2-\ff dp-\ff 2 q)$,
$$\|u\|_\infty+ \|\nn u\|_\infty  \le \ff 1 3,\ \     \| u\|_{\tt W_\infty^{\theta,\infty}(T)}+ \|\nn^2 u_t\|_{\tt L^{p_1'}_{q_1'}(T)}<\infty,$$
where for some $0\le h\in C_0^\infty(\R^d)$ with $h|_{B(0,1)}=1,$
$$\| u\|_{\tt W_\infty^{\theta,\infty}(T)} :=  \sup_{(x,t)\in\R^d\times [0,T]} \|u_t h(x+\cdot) \|_{W^{\theta,\infty}}.$$
By the Sobolev embedding theorem,       $\| u\|_{\tt W_\infty^{\theta,\infty}(T)}<\infty$ for some $\theta>1$ implies that $\nn u_t$ is H\"older continuous uniformly in $t\in [0,T].$

Let $$\Theta_t(x):= x+ u_t(x),\ \ \ (t,x)\in [s,T]\times \R^d.$$  Then $\Theta_t$ is diffeomorphism uniformly in $t\in [s,T]$. By \eqref{PL0} and    It\^o's formula  \cite[Lemma 1.2.3(3)]{RW24},
$X_{s,t}^x$ solves \eqref{S0} if and only if
$ Y_{s,t}^x:= \Theta_t(X_{s,t}^x)$ solves the SDE
 \beq\label{YC}\d Y_{s,t}^x = \bar{b}_t(Y_{s,t}^x)+ \bar \si_t(Y_{s,t}^x)\d W_t,\ \ Y_{s,s}^x= x+ u_{s}(x),\ t\in [s,T],\end{equation}
where
\beq\label{XC} \bar{b}_t:= \left(\ll u_t + b_t^{(1)}+\sum_{i=2}^{\ell'} (\nn \Theta_t)  b_t^{0,i}\right)\circ\Theta_t^{-1},\ \ \bar \si_t:= \big\{(\nn \Theta_t)\si_t\big\}\circ\Theta_t^{-1}.\end{equation}
By $(C)$ and the properties of $u$ mentioned above,  we see that   $(C)$ with  $\ell'=n$  holds for $(\bar b, \bar \si)$ in place of $(b,\si)$.   So, the  assumption on $\ell'=n$ implies that   \eqref{YC} is weakly (also strongly
when $(A_2)$ holds) well-posed,  and  so is \eqref{S0} for $X_{s,t}^x= \Theta_t^{-1}(Y_{s,t}^x)$.   \end{proof}

  \subsection{Regularities}

When $(C)$ and $(A_2)$ hold  with $\ell'=1$, the  moment estimates, log-Harnack inequality and Bismut formula have been derived, see \cite[Theorems 1.3.1, 1.4.2, 1.5.1]{RW24}. The next result extend these
to the case $\ell'\ge 2$.
Moreover,   we formulate these estimates with   explicit  dependence on $\|  b^{0,i}\|_{\tt L_{q_i'}^{p_i'}(s,t)},$   which is crucial
in the proof of Theorem \ref{T02}(2).

\beg{prp} \label{L*1} Assume $(C)$ and $(A_2)$.   Then   for any $q\in [1,\infty)$, there exist   constants $c,l\in [2,\infty)$
depending only on $d,K, T, p_i',q_i'$ and $a$, such that the following assertions hold for any $(s,x)\in [0,T)\times \R^d,$ and any $t\in [s,T]$.
\beg{enumerate}
\item[$(1)$] For any $(p',q')\in \scr K$ and $g\in \tt L_{q'}^{p'}(0,t)$,
\beq\label{PG2'} \E\big[\e^{\int_0^t g^2(X_{s}^x)\d s}\big]  \le c  \exp\bigg[c \sum_{i=1}^{\ell'} \|  b^{0,i}\|^l_{\tt L_{q_i'}^{p_i'}(t)}+ c  \|g\|_{\tt L_{q'}^{p'}(t)}^l\bigg].\end{equation}
 \item[$(2)$] There holds
 \beq\label{MM} \E\bigg[\sup_{t\in [s,T]}|X_{s,t}^{x}|^q\bigg]\le c\exp\bigg[c\sum_{i=1}^{\ell'}\|  b^{0,i}\|^2_{\tt L_{q_i'}^{p_i'}(s,T)}\bigg] (1+|x|^q).\end{equation}
 \item[$(3)$] For any $v\in\R^d,$
 $$\nn_v X_{s,t}^x:= \lim_{\vv\downarrow 0} \ff{X_{s,t}^{x+\vv v}-X_{s,t}^x}\vv$$
 exists in $L^q(\OO\to\R^d,\P)$, and
 \beq\label{DR} \E\bigg[\sup_{t\in [s,T]} |\nn_v X_{s,t}^x|^q\bigg] \le c|v|^q\exp\bigg[c\sum_{i=1}^{\ell'}\|  b^{0,i}\|^l_{\tt L_{q_i'}^{p_i'}(s,T)}\bigg].\end{equation}
\item[$(4)$] The following log-Harnack inequality holds  for $f\in \B_b^+(\R^d),$ $0\le s<t\le T$ and $x,y\in\R^d$:
\beq\label{KK1}  P_{s,t} \log f(x)\le \log P_{s,t} f(y)+\ff {c|x-y|^2} {t-s}\exp\bigg[c\sum_{i=1}^{\ell'}\|  b^{0,i}\|^l_{\tt L_{q_i'}^{p_i'}(s,t)}\bigg], \end{equation} where $\B_b^+(\R^d)$ is the set of all strictly positive bounded measurable functions on $\R^d$.
\item[$(5)$] For any   $v \in\R^d,$ $\bb\in C^1([s,t])$ with $\bb_s=0,\bb_t=1$,
\beq\label{BSM} \nn_v P_{s,t}f(x)=\E\bigg[f(X_{s,t}^x)  \int_s^t \bb_r'\big\<\{\si_r^*a_r^{-1}\}(X_{s,r}^x)\nn_v X_{s,r}^x,\ \d W_r\big\>\bigg].\end{equation}
\end{enumerate}
\end{prp}

\beg{proof} By \eqref{DR} for $q=2$, we obtain
$$|\nn P_{s,t}f|^2\le c  (P_{s,t}|\nn f|^2) \exp\bigg[c\sum_{i=1}^{\ell'}\|  b^{0,i}\|^l_{\tt L_{q_i'}^{p_i'}(s,t)}\bigg],$$
which  implies
  the log-Harnack inequality \eqref{KK1} for some possibly different constant $c\in (0,\infty)$, see the proof of \cite[Theorem 1.5.1]{RW24} or \cite[Proof of (2.18)]{FYW2} for details.
Hence, we only need to prove \eqref{PG2'}, \eqref{MM},  \eqref{DR} and \eqref{BSM}.

Without loss of generality, we only consider
$s=0$, and denote $X_t^x=X_{0,t}^x,P_{t}=P_{0,t}$. All constants below depend only on  $d, p_i',q_i',K, T$ and $a.$

(a) Let $\ell'=1$. As indicated above that  the well-posedness, the existence of $\nn_v X_{s,t}^x$ and the Bismut formula   \eqref{BSM} are already known.
As explained on the equation \eqref{PL0},   for any $\ll\in (0,\infty)$ the PDE
\beq\label{PL} (\pp_s +L_s-\ll)u_s=-  b_s^{0,1},\ \ \ s\in [0,t],\ u_t=0\end{equation}
has a unique solution $u: [0,t]\times \R^d\to\R^d$ such that
$$\|u\|_\infty+ \|\nn u\|_\infty+\|u\|_{\tt W_{\infty}^{\theta,\infty}(t)} + \|\nn^2 u\|_{\tt L_{q_1'}^{p_1'}(t)}<\infty, $$
 where $\theta\in (1,2-\ff d p-\ff 2 q).$

  To estimate the upper bound using $\|  b^{0,1}\|_{\tt L_{q_1'}^{p_1'}(t)},$ let
\beg{align*}&\bar L_s:= \ff 1 2{\rm tr}\{ a_s\nn^2\}+    b^{(1)}_s \cdot\nn,\ \ \
  \bar u_s:= \ff{u_s}{1+\|  b^{0,1}_s\|_{\tt L_{q_1'}^{p_1'}(t)}},\\
&f_s= \ff{  b^{0,1}_s}{1+\|  b^{0,1}\|_{\tt L_{q_1'}^{p_1'}(t)}} +\ff{  b^{0,1}_s} {1+\|  b^{0,1}\|_{\tt L_{q_1'}^{p_1'}(t)}}\cdot\nn u_s, \ \ \ s\in [0,t].\end{align*}
Then
 \beg{align*} &(\pp_s +\bar L_s-\ll)\bar u_s  =-f_s,\ \ \ s\in [0,t],\ \bar u_t=0,\\
&\|f\|_{\tt L_{q_1'}^{p_1'}(t)}\le  1+\|\nn u\|_\infty.\end{align*}
Combining this with  \cite[Lemma 1.2.2]{RW24}, we find   constants $c_1, \ll_0\geq1 $ such that when $\ll\ge \ll_0,$
\beg{align*} &\ff{\|u\|_\infty}{1+\|  b^{0,1}\|_{\tt L_{q_1'}^{p_1'}(t)}}=\|\bar u\|_\infty\le c_1 \ll^{-1}\|f\|_{\tt L_{q_1'}^{p_1'}(t)}\le  c_1 \ll^{-1} \big(1+\|\nn u\|_\infty\big),\\
& \ff{\|\nn   u\|_\infty}{1+\|  b^{0,1}\|_{\tt L_{q_1'}^{p_1'}(t)}}=\|\nn\bar u\|_\infty\le c_1 \ll^{-\ff 1 2}\|f\|_{\tt L_{q_1'}^{p_1'}(t)}\le  c_1 \ll^{-\ff 1 2} \big(1+\|\nn u\|_\infty\big),\\
&\ff{\|\nn^2  u\|_{\tt L_{q_1'}^{p_1'}(t)}}{1+\|  b^{0,1}\|_{\tt L_{q_1'}^{p_1'}(t)}}= \|\nn^2\bar u\|_{\tt L_{q_1'}^{p_1'}(t)} \le c_1 \|f\|_{\tt L_{q_1'}^{p_1'}(t)}\le c_1 \Big(1+\|\nn u\|_\infty\Big).\end{align*}
Taking
\beq\label{PPL}\ll:=\ll_0\lor \Big(4 c_1 \big(1+\|  b^{0,1}\|_{\tt L_{q_1'}^{p_1'}(t)}\big)\Big)^2,\end{equation}
we derive
\beq\label{67} \|u\|_\infty\lor\|\nn u\|_\infty\le \ff 1 3,\ \ \|\nn^2 u\|_{\tt L_{q_1'}^{p_1'}(t)}\le 2 c_1
\Big(1+\|  b^{0,1}\|_{\tt L_{q_1'}^{p_1'}(t)}\Big).\end{equation}

Let
$$\Theta_s(x):=x+u_s(x),\ \ Y_{s}^x:= \Theta_s(X_{s}^x),\ \ s\in [0,t], \ x\in\R^d.$$
By \eqref{67} we have
\beq\label{76} \beg{split}&\|\nn \Theta_s\|_\infty+\|(\nn \Theta_s)^{-1}\|_\infty\le 2,\\
&\ff 1 2 |X_s^x-X_{s}^y|\le |Y_{s}^x-Y_{s}^y|\le 2 |X_{s}^x-X_{s}^y|,\ \ s\in [0,t],\ x,y\in\R^d.\end{split}\end{equation}
Similarly to \eqref{YC}, we have
\beq\label{YX}\d Y_{s}^x = \bar{b}_s(Y_{s}^x)+ \bar \si_s(Y_{s}^x)\d W_s,\ \ Y_{0}^x= x+ u_{0}(x),\ s\in [0,t],\end{equation}
where
$\bar{b}$ and $\bar \si$ are in \eqref{XC} with $\ell'=1$.  By
 \eqref{67}, we find a  constant $c_5\in (0,\infty)$ such that
\beq\label{770}\beg{split}&\|\nn \bar{b}_s\|_\infty:=\sup_{x\neq y}\frac{|\bar{b}_s(x)-\bar{b}_s(y)|}{|x-y|}\le c_5(\ll+ 1),\\
& \|\nn \bar \si_s\|^2\le c_5\big(\|\nn \si_s\|^2+  \|\nn^2 u_s\|^2\big)\circ\Theta_s^{-1},\ \  \ s\in [0,t].\end{split}\end{equation}
By   Krylov's and Khasminski's estimates, see \cite{YZ} or   \cite[Theorem 1.2.3 (2), Theorem 1.2.4]{RW24}
for $\d A_s=\d s$, we find  constants $c_2, l\ge 2$ such that
$$\E\bigg[\e^{\int_0^t g^2(Y_{s}^x)\d s}\bigg]\le c_2\exp\bigg[c_2 \|g\|_{\tt L_{q'}^{p'}(t)}^l\bigg],\ \ g\in \tt L_{q'}^{p'}(0,t).$$
Combining this together with \eqref{PPL}, \eqref{67} and
$$\E\bigg[\e^{\int_0^t g^2(X_{s}^x)\d s}\bigg]=\E\bigg[\e^{\int_0^t (g\circ\Theta_s^{-1})^2(Y_{s}^x)\d s}\bigg],$$
we prove \eqref{PG2'} for $\ell'=1$.

Next, by It\^o's formula  and the maximal functional inequality \cite[Theorem 2.1]{XXZZ}, for any $q\ge 2$, we find a constant $c_6\in (0,\infty)$ such that
\beg{align*}&\d |Y_{s}^x-Y_{s}^y|^{2q}\le  c_6 |Y_{s}^x-Y_{s}^y|^{2q} \big\{1+\ll + \scr M \|\nn \bar \si_s\|^2 (Y_{s}^x)+
\scr M \|\nn \bar \si_s\|^2 (Y_{s}^y) \big\}\d s+\d M_s,\\
&\d |Y_{s}^x|^{2q}\le  c_6 \big(1+|Y_{s}^x|^{2q} \big) (1+\ll)\d s+\d N_s \ \ \ s\in [0,t],\end{align*}
for some martingales $M_s $  and $N_s$.  Thus, by the  stochastic Gronwall lemma and  maximal functional inequality (see  \cite[Lemma 1.3.3, Lemma 1.3.4]{RW24} or \cite{XXZZ}),   and applying \eqref{PG2'}, \eqref{67} and \eqref{770}, we find a constant $c_7\in (0,\infty)$ such that
\beg{align*}&\E\bigg[\sup_{s\in [0,t]} |Y_{s}^x-Y_{s}^y|^q\bigg]\le c_7 \exp\bigg[c_7\ll + c_7\|b^{0,1}\|_{\tt L_{q_1'}^{p_1'}(t)}^l\bigg]\big|\Theta_0(x)-\Theta_0(y)\big|^q,\\
&\E\bigg[\sup_{s\in [0,t]} |Y_{s}^x|^q\bigg]\le c_7 \e^{c_7\ll}\Big(1+\big|\Theta_0(x)\big|^q\Big),  \  \ x,y\in\R^d.\end{align*}
This together with \eqref{PPL}, \eqref{76} and $l\ge 2$ yields that for some constant $c_8\in (0,\infty)$,
\beg{align*}& \E\bigg[\sup_{s\in [0,t]} |X_{s}^x-X_{s}^y|^q\bigg]
\le c_8 \exp\bigg[c_8 \|b^{0,1}\|_{\tt L_{q_1'}^{p_1'}(t)}^l\bigg]|x-y|^q,\\
&\E\bigg[\sup_{s\in [0,t]} |X_{s}^x|^q\bigg]
\le c_8 \exp\bigg[c_8 \|b^{0,1}\|_{\tt L_{q_1'}^{p_1'}(t)}^2\bigg]\big(1+|x|^q\big) \  \ x,y\in\R^d.\end{align*}
 The second estimate implies \eqref{MM}, while the first estimate together  with the definition of $\nn_v X_{s,t}^x$ implies  \eqref{DR},  for $c=c_8$.

 (b) Assume that the assertions hold for $\ell'=n$ for some $n\in \mathbb N$.
 We consider the case
 for $\ell'=n+1$.

 Let  $u_s$ and $\Theta_s$ be constructed above for $\ll$ satisfying  \eqref{PPL}.  By It\^o's formula,
   $Y_s^x:= \Theta_s(X_s^x)$ solves the SDE \eqref{YX}, where
   as explained in the proof of Proposition \ref{P01} that the coefficients of this SDE satisfy $(C)$ for $\ell'=n$.   So, by the induction assumption,
 all assertions hold for $Y_s^x$ in place of $X_s^x$, which together with $Y_s^x=\Theta_s(X_s^x)$, \eqref{67} and \eqref{76}, imply  estimates \eqref{MM}-\eqref{DR}
 for some constant $c\in (0,\infty)$. Moreover, by
   chain rule and
$$ P_t f(x)=\bar P_t(f\circ \Theta_t)(\Theta_0(x))  :=\E[(f\circ \Theta_t)(Y_t^x)\big],$$
  \eqref{BSM} follows from the corresponding  formula for $\bar P_t,$
   see \cite[page 32]{RW24} or \cite[page 1876]{FYW3} for details.
\end{proof}

\subsection{Local hyperbound  estimates}

We  first consider the local hyperbound  on the diffusion semigroup
$$\hat P_{s,t}f(x):=\E[f(\hat X_{s,t}^x)],\ \ \ 0\le s\le t<\infty,\ f\in \B_b(\R^d)$$
 associated with the SDE
\beq\label{HSDE} \d \hat X_{s,t}^x= \hat b_t(\hat X_{s,t}^x)\d t + \si_t(\hat X_{s,t}^x)\d W_t,\ \
t\ge s,\ X_{s,s}^x=x, \end{equation}
where the noise coefficient $\si$ and drift  $ \hat b$  satisfy  the following assumption.

 \beg{enumerate} \item[$(A_1')$] For any $T\in (0,\infty)$, $a:=\si\si^*$ satisfies the corresponding condition in $(A_1)$ for some constants    $K\in (0,\infty)$ and $\aa\in (0,1]$, and moreover
$$   \|\hat b_t(0)\|\le K,\ \
 |\hat b_t(x)-\hat b_t(y)|\leq K(1+|x-y|),\ \ t\in [0,T],\ x,y\in \R^d. $$
 \end{enumerate}

 It is well known that under $(A_1')$, for any $s\in [0,\infty)$ and $x\in \R^d$, the SDE
 \eqref{HSDE}
is weakly well-posed, see for instance \cite{BK}.

\beg{lem}\label{L2'} Assume $(A_1')$.
 Then for any $T\in (0,\infty)$ there exists a constant
 $c\in (0,\infty)$ depending only on $(d,T,K,\aa)$ such that for any $1\le p_1\le p_2\le \infty$ and $ 0\le s<t\le T,$
 \beq\label{ES1} \|\hat P_{s,t}\|_{\tt L^{p_1}\to \tt L^{p_2}}:=\sup_{\|f\|_{\tt L^{p_1}}\le 1}
 \|\hat P_{s,t}f\|_{\tt L^{p_2}}
  \le c (t-s)^{-\ff{d(p_2-p_1)}{2p_1p_2}},\end{equation}
  \beq\label{ES2} \|\nn\hat P_{s,t}\|_{\tt L^{p_1}\to \tt L^{p_2}}:=\sup_{\|f\|_{\tt L^{p_1}}\le 1}
 \|\nn\hat P_{s,t}f\|_{\tt L^{p_2}}
  \le c (t-s)^{-\ff 1 2-\ff{d(p_2-p_1)}{2p_1p_2}}.\end{equation}
\end{lem}

\beg{proof}  The desired estimates for  $L^{p_1}$-$L^{p_2}$  in  place of $\tt L^{p_1}$-$\tt L^{p_2}$ are well known.  The proof for the present estimates is based on a localization argument
as in \cite{FYW23JDE}.  All constants below depend only on $(d,T,K,\aa)$ in $(A_1')$.

By  \cite[Theorem 1.2 (I)-(II)]{MPZ}, the condition  $(A_1')$ implies that
 $\hat P_{s,t}$ has density $\hat{p}_{s,t}(x,y)$ with respect  to the Lebesgue measure such that
for some  constants $c_0,\kk\in (0,\infty)$,
\beq\label{G1} \hat p_{s,t}(x,y)\le c_0 p_{t-s}^\kk(\psi_{s,t}(x)-y),\end{equation}
\beq\label{G2} |\nn \hat p_{s,t}(\cdot,y)|(x)\le c_0 (t-s)^{-\ff 1 2} p_{t-s}^\kk(\psi_{s,t}(x)-y)\end{equation}
hold for all $0\le s<t\le T$ and $x,y\in\R^d$, where
$$p_t^\kk(z):= (\kk\pi t)^{-\ff d 2}\e^{-\ff{|z|^2}{\kk t}},\ \ t>0,\ z\in\R^d,$$
and $\{\psi_{s,t}\}_{0\le s\le t\le T} $  is a family of diffeomorphisms on $\R^d$ satisfying
\beq\label{3} \sup_{0\le s\leq t\le T}\big\{\|\nn \psi_{s,t}\|_\infty+\|\nn \psi^{-1}_{s,t}\|_\infty
\big\}\le \dd \end{equation} for some constant $\dd\in (0,\infty)$.

Let
$$P_t^\kk f(x):=\int_{\R^d}p_t^\kk(x-y)f(y)\d y, \ \ f\in \B_b(\R^d), x\in\R^d,t\geq 0.$$
It is classical that for some constant $c(\kk,d)\in (0,\infty)$
\beq\label{UP} \|\nn^i P_t^\kk\|_{L^{p_1}\to L^{p_2}}\le c(\kk,d) t^{-\ff i 2-\ff {d(p_2-p_1)}{2p_1p_2}},\ \ t>0,\ i=0,1,\end{equation}
where   $\nn^0 $ is the identity operator. For any $n\in \Z_+$, let
$${\bf B}_n:=\bigg\{v\in \Z^d:\ |v|_1:=\sum_{i=1}^d |v_i|=n\bigg\}.$$
Then for any $z\in\R^d$,
\beq\label{UP2} 1 \le \sum_{n=0}^\infty\sum_{v\in {\bf B}_n} 1_{B(z+v,d)},\end{equation}
where $B(y,d):=\{x\in\R^d:\ |x-y|\le d\}$.
Moreover, by \eqref{3}, we find a constant $c_1>1$ such that
\beg{align*}&|\psi_{s,t}(x)-y|^2\ge c_1^{-1} n^2 -c_1,\\
 &\   x\in B(\psi_{s,t}^{-1}(z),1),\ y\in \cup_{v\in {\bf B}_n} B(z+v,d),\ z\in\R^d,\ \ 0\le s\le t\le T,  n\in\Z_+.\end{align*}
So, there exists a constant $c_2>1$ such that
\beg{align*}&p_t^\kk(\psi_{s,t}(x)-y)\le c_2 \e^{-c_2^{-1}n^2} p_t^{2\kk}(\psi_{s,t}(x)-y),\\
& \ x\in  B(\psi_{s,t}^{-1}(z),1),\ y\in \cup_{v\in {\bf B}_n} B(z+v,d),\ 0\le s\le t\le T, \ n\in\Z_+.\end{align*}
Combining this with \eqref{G1}, \eqref{3} and \eqref{UP}, we find constants $ c_3,c_4\in (0,\infty)$ such that
  \beg{align*} &\big\|1_{B(\psi_{s,t}^{-1}(z),1)}\hat P_{s,t}f\|_{L^{p_2}}\le c_0 \sup_{\|g\|_{L^{\ff{p_2}{p_2-1}}}\le 1} \int_{\R^d\times\R^d}|g1_{B(\psi_{s,t}^{-1}(z),1)}|(x)p_{t-s}^\kk(\psi_{s,t}(x)-y)|f|(y)\d x\d y\\
 & \le c_0 \sup_{\|g\|_{L^{\ff{p_2}{p_2-1}}}\le 1}\sum_{n=0}^\infty \sum_{v\in {\bf B}_n}
 \int_{\R^d\times\R^d}|g1_{B(\psi_{s,t}^{-1}(z),1)}|(x)p_{t-s}^\kk(\psi_{s,t}(x)-y)(|f|1_{B(z+v,d)})(y)\d x\d y \\
 &\le c_3  \sum_{n=0}^\infty n^{d-1}\e^{-c_2^{-1} n^2} \sup_{\|g\|_{L^{\ff{p_2}{p_2-1}}}\le 1}\sup_{v\in \Z^d}
 \int_{\R^d}|g| (x)P_{t-s}^{2\kk} (|f|1_{B(z+v,d)})(\psi_{s,t}(x))\d x \\
 &\le c_4 \sup_{v\in\Z^d}\|P_{t-s}^{2\kk} (|f|1_{B(z+v,d)})\|_{L^{p_2}}\le c_5 (t-s)^{-\ff{d(p_2-p_1)}{2p_1p_2}} \|f\|_{\tt L^{p_1}},\ \ z\in\R^d, f\in \B_b(\R^d).\end{align*}
 Taking supremum over $z\in\R^d$ we prove
  \eqref{ES1}   for some constant $c\in (0,\infty)$. The estimate \eqref{ES2} can be proved in the same way
 by using \eqref{G2} in place of \eqref{G1}.
\end{proof}

We are now ready to prove   the following result.

\beg{prp}\label{P01'} Assume $(C)$ and there exist constants $\tilde{K}>0$ and $\alpha\in(0,1]$ such that \begin{align}\label{aap}|a_t(x)-a_t(y)|\le \tilde{K}|x-y|^\aa,\ \ t\in[0,T], x,y\in\R^d.
 \end{align}Then for any $1< p_1\le p_2\le \infty$, there exists a constant $c\in (0,\infty)$ depending only on $T,d,K,\tilde{K}$ and $\alpha$,
such that
\beq\label{ES1''} \| P_{s,t} \|_{\tt L^{p_1}\to\tt L^{p_2}}\le c(t-s)^{-\ff{d(p_2-p_1)}{2p_1p_2}},
\ \ 0\le s< t\le T,\end{equation}
\beq\label{ES2''} \|\nn P_{s,t} \|_{\tt L^{p_1}\to\tt L^{p_2}}\le c(t-s)^{-\ff 1 2-\ff{d(p_2-p_1)}{2p_1p_2}},
\ \ 0\le s< t\le T.\end{equation} When $ b^{0,i}=0$ for $1\le i\le \ell'$, these estimates also hold for $p_1=1.$
\end{prp}

\beg{proof}    Let $\hat P_{s,t}$ and $\hat X_{s,t}^x$ be in Lemma \ref{L2'} for
$$\hat b_t(x):= b_t^{(1)}(x).$$
By $(C)$ and \eqref{aap}, the condition $(A_1')$ in Lemma \ref{L2'} holds, so that \eqref{ES1} and \eqref{ES2} hold for some constant $c\in (0,\infty)$.
When  $ b^{0,i}=0$ for $1\le i\le \ell'$, we have $\hat P_{s,t}= P_{s,t},$ so that \eqref{ES1''} and \eqref{ES2''} hold for any $1\le p_1\le p_2\le\infty.$

 In general, by $(C)$, for any $r\in [s,T]$,
$$\xi_r^i:= \big\{\si_r^*a_r^{-1}(\hat X_{s,r})\big\}  b_r^{0,i} (\hat X_{s,r}),\ \ 1\le i\le \ell'
$$ satisfies
$$|\xi_r^i|\le c_1  |b_r^{0,i}|(\hat X_{s,r}),\ \ \|b^{0,i}\|_{\tt L_{q_i'}^{p_i'}(T)}\le K $$
for some constant $c_1\in (0,\infty).$
Then by Krylov's and Khasminskii's estimates, see e.g. \cite[Theorem 1.2.3 and Theorem 1.2.4]{RW24}, and H\"{o}lder's inequality,
we have
\begin{align*}
&\E\e^{ q\int_s^t |\sum_{i=1}^{\ell'}\xi_r^i|^2\d r}\leq \E\prod_{i=1}^{\ell'}\e^{ q\ell'\int_s^t |\xi_r^i|^2\d r}\leq \prod_{i=1}^{\ell'}\left(\E\e^{ q(\ell')^2\int_s^t |\xi_r^i|^2\d r}\right)^{\frac{1}{\ell'}}\leq c_1(q),\ \ t\in [s,T],q>0.
\end{align*}
This means that
$$R_t:=\e^{\int_s^t\<\sum_{i=1}^{\ell'}\xi_r^i,\d W_r\>-\ff 1 2 \int_s^t |\sum_{i=1}^{\ell'}\xi_r^i|^2\d r},\ \ t\in [s,T]$$
is a martingale, and there exists a constant $c_2>1$ such that for $\tt p_1:=\ss{p_1},$
\begin{align*}&\Big(\E\Big[R_t^{\ff {\tt p_1}{\tt p_1 -1}}\Big]
\Big)^{\ff{\tt p_1-1}{\tt p_1 }}\leq \Big(\E\Big[\e^{\ff {\tt p_1^2+\tt p_1}{(\tt p_1 -1)^2}\int_s^t |\sum_{i=1}^{\ell'}\xi_r^i|^2\d r
}\Big]
\Big)^{\ff{\tt p_1-1}{2\tt p_1 }}
\le c_2,\ \ t\in [s,T].
\end{align*}
So, by   Girsanov's theorem and H\"older's inequality, we obtain
$$|  P_{s,t}  f|=\bigg|\E\Big[f(\hat X_{s,t})R_t\Big] \bigg|
\le c_2 \big(\hat P_{s,t}|f|^{\tt p_1}\big)^{1/\tt p_1}.$$
This together with \eqref{ES1} yields that for some constant $c_3\in (0,\infty)$
\beg{align}\label{Hyp12}\nonumber&\|  P_{s,t}  f\|_{\tt L^{p_2}}
\le c_2 \big\|\hat P_{s,t}|f|^{\tt p_1}\big\|_{\tt L^{p_2/\tt p_1}}^{1/\tt p_1}\\
&\le c_2 \|f\|_{\tt L^{p_1}} \|\hat P_{s,t}\|_{\tt L^{\tt p_1}\to \tt L^{p_2/\tt p_1}}^{1/\tt p_1}\le c_3 \|f\|_{\tt L^{p_1}} (t-s)^{-\ff{d(p_2-p_1)}{2p_1p_2}},\ \ 0\le s<t\le T.\end{align}
Thus, \eqref{ES1''} holds for some constant $c\in (0,\infty).$

By \eqref{DR} and \eqref{BSM} in Proposition \ref{L*1}, we find a constant $c_4\in (0,\infty)$
such that
\begin{align}\label{nab}|\nn  P_{s,t}  f|\le c_4 (t-s)^{-\ff 1 2} \big(P_{s,t}  |f|^{\tt p_1}\big)^{1/\tt p_1},\ \ 0\le s<t\le T.
\end{align}
 Combining this with \eqref{ES1'} for $(\tt p_1,  p_2/\tt p_1)$ in place of $(p_1,p_2)$, we find a constant $c_5\in (0,\infty)$ such that for any $0\le s<t\le T$,
 $$\|\nabla  P_{s,t}\|_{\tt L^{p_1}\to \tt L^{p_2}}\le c_4 (t-s)^{-\ff 1 2 } \|  P_{s,t}\|_{\tt L^{\tt p_1 }\to \tt L^{p_2/\tt p_1}}^{1/\tt p_1 }\le c_5 (t-s)^{-\ff 1 2- \ff{d(p_2-p_1)}{2p_1p_2}}.$$
 Hence,     \eqref{ES2''} holds for   $c=c_5$.
\end{proof}

\subsection{Comparing two singular SDEs}

Next,  we  estimate   the distance of solutions to   different SDEs.
For $b^{0,i}$ in $(C)$, and let
\beq\label{CV0}\tt b^{0,j}\in \tt L_{\tt q_j}^{\tt p_j}(T)\ \text{ for\ some\ }  \{(\tt p_j,\tt q_j)\}_{1\le j\le \tt \ell}\subset \scr K.\end{equation}  We  denote
$$b^{0}:= \sum_{i=1}^{\ell'} b^{0,i},\ \ \  \tt b^{0}=\sum_{j=1}^{\tt\ell} \tt b^{0,j}.$$
Consider   the SDEs \eqref{S0} and
\beq\label{S0'}\d   \tt X_{s,t}^y= \big(\tt b_t^{0}+   b_t^{(1)}\big)(X_{s,t}^y)\d t +\si_t(X_{s,t}^y)\d W_t,\ \ t\in [s,T],\ \tt X_{s,s}^y=y.\end{equation}
We have the following result.

\beg{prp} \label{L*2} Assume $(C)$,   $(A_2)$ and $\eqref{CV0}$. Let $X_{s,t}^x$ and $\tt X_{s,t}^y$ solve
\eqref{S0} and \eqref{S0'} respectively.
\beg{enumerate}
\item[$(1)$]
For any   $q\in [1,\infty)$, we find  constants $c,l\ge 2$ depending only on $d, p_i', q_i',\tt p_j, \tt q_j, K, T$ and $\zeta$ in $\eqref{conti}$, such that
 \beq\label{KK2} \begin{split}&\E\bigg[\sup_{s\in[r, t]}|X_{r,s}^{x}-  \tt X_{r,s}^{y}|^q\bigg] \le \exp\bigg[c+c  \sum_{i=1}^{\ell'} \|b^{0,i}\|_{\tt L_{q_i'}^{p_i'}(r,t)}^l  \bigg] \\
&\qquad \times
  \bigg(|x-y| + \int_r^t   \| b_s^{0}   -   \tt b_s^{0}\|_{\tt L^\infty}\d s \bigg)^q,\ \
  0\le r<t\le T,\ x,y\in\R^d.\end{split}\end{equation}
 \item[$(2)$] For any $f\in \B_b(\R^d)$,
\begin{align}\label{DH1'}  \tt P_{r,t} f=  P_{r,t}  f+\int_r^t \tt P_{r,s}  \<\tt b_s^{0}-  b_s^{0}, \nn  P_{s,t}  f\>\d s, \ \ 0\le r\le t\le T.
\end{align}\end{enumerate}
  \end{prp}
\beg{proof} Without loss of generality, we simply let $r=0.$

   (1) Denote $X_{0,t}^{x}=X_t^x,   \tilde{X}_{0,t}^{y}=   \tilde{X}_t^{y}$ for $t\in [0,T]$ and $x,y\in \R^d$.   Similarly to   step (b) in the proof of Proposition \ref{L*1} by inducing in $\ell'$, we only need to prove  the desired assertion  for $\ell'=1.$

Let $\ell'=1$, and simply denote
$$k_t:=1+  \|b^{0,1}\|_{\tt L_{q_1'}^{p_1'}(t)}^l+
\sum_{j=1}^{\tt\ell}\|\tt b^{0,j}\|_{\tt L_{\tt q_j}^{\tt p_j}(t)}^l$$
for $  t\in (0,T],$ where $l\ge 2$ is the constant in Proposition \ref{L*1}.
All constants $c_j \ge 2$ below   depend only on $d, p_1', q_1', \tt p_j, \tt q_j,  K, T$ and $\zeta$.

Let $\ll$ be in \eqref{PPL} such that \eqref{67} holds for  $u$ solving \eqref{PL}.
For fixed $t\in (0,T],$ let
 \beq\label{7} Y_s^{x}=  X_s^{x} + u_s( X_s^{x}),\ \ \ \   \tt Y_s^{y}=  \tt X_s^{y}+u_s(\tt X_s^{y}),\ \ \ s\in [0,t], x,y\in \R^d.\end{equation}
 By \eqref{76} we have
 \beq\label{77} |X_s^{x}-  \tt X_s^{y}|\le 2 |Y_s^{x}-  \tt Y_s^{y}|,\ \ \ s\in [0,t], x,y\in \R^d.\end{equation}
By It\^o's formula,  \eqref{YX} holds and
 $$\d  \tt Y_{s}^{y}= \big\{ \bar{b}_s+ \big\{(\nabla\Theta_s)(\tt b_s^{0}-   b_s^{0})\big\}\circ\Theta_s^{-1}\big\}(\tt Y_{s}^{y})\d s
+\bar\si_s(\tt Y_{s}^{y})\d W_s,\ \ s\in [0,t].$$
Combining this with \eqref{YX}, \eqref{77}, \cite[Lemma 2.1]{XXZZ}  and applying It\^o's formula,  for fixed $q\in [1,\infty)$,   we find a constant $c_1\in (0,\infty)$
and a martingale $M_t$ such that
\beg{align*}  \d |Y_{s}^x-\tt Y_{s}^y|^{q+1} \le &\, c_1 |Y_{s}^{x}-  \tt Y_{s}^{y}|^{q+1} \big(1+\ll + \scr M \|\nn \bar{\si}\|^2 (Y_{s}^{x})+\scr M \|\nn \bar{\si}\|^2 (\tt Y_{s}^{y})\big)\d s\\
&+c_1 |Y_{s}^{x}- \tt Y_{s}^{y}|^{q}\|b_s^{0}-  \tt b_s^{0}\|_\infty  \d s+\d M_s,\ \ \ s\in [0,t].\end{align*}
By \eqref{PG2'} for $Y_s^x$ and $\tt Y_s^y$ in place of $X_s$, the stochastic Gronwall lemma, the maximal function inequality and Khasminski's estimate as explained above, see for instance \cite[Lemma 1.3.3]{RW24} and \cite[Theorems 1.2.3, 1.2.4]{RW24}, we find   constants $c_2,c_3\in (0,\infty)$ such that
\beg{align*} &\bigg(\E\bigg[\sup_{s\in [0,t]} |Y_{s}^{x}- \tt Y_{s}^{y}|^q\bigg]\bigg)^{\ff{q+1}q}
 \le   \e^{c_2k_t}
 \bigg(|x-y|^{q+1}+ \E\int_0^t |Y_{s}^x-\tt Y_{s}^y|^{q}\| b_s^{0}- \tt b_s^{0}\|_\infty\d s\bigg)\\
&\le   \e^{c_2k_t}
 \bigg(|x-y|^{q+1}+
  \E\bigg[\sup_{s\in [0,t]}|Y_{s}^{x}-  \tt Y_{s}^{y}|^{q}\bigg]\int_0^t\|b_s^{0}- \tt b_s^{0}\|_\infty\d s \bigg) \\
&\le
     \ff 1 2 \bigg(\E\bigg[\sup_{s\in [0,t]} |Y_{s}^{x}- \tt Y_{s}^{y}|^q\bigg]\bigg)^{\ff{q+1}q}+    \e^{c_3k_t}
   \bigg(|x-y| +\int_0^t\|b_s^{0}-  \tt b_s^{0}\|_\infty\d s \bigg)^{q+1}.
\end{align*}
  Combining this with   \eqref{77},
   we obtain \eqref{KK2}
for some constant  $c \ge 2$.

(2) By \eqref{nab} and an approximation argument, it suffices to prove the desired assertion for  $f\in C_0^\infty(\R^d)$.  We first prove the Kolmogorov backward equation
\beq\label{BC} \pp_s P_{s,t}  f= - L_s P_{s,t} f,\ \ 0\le s\le t\le  T,\  f\in C_0^\infty(\R^d).\end{equation}
Let  $\ell'=1$, then
$$L_s =  \ff 1 2 {\rm tr}(a_s\nn^2)+ (b^{(1)}_s+b^{0,1}_s)\cdot\nn,\ \ s\in [0,T].$$
By It\^o's formula we obtain the forward   Kolmogorov equation
\beq\label{FW} \pp_t  P_{s,t}   f =   P_{s,t}  L_t  f,\ \ s\le t\le T. \end{equation}
  By $(C)$, for any  $f\in C_0^\infty(\R^d)$, we have $\|L_s   f\|_{\tt L_{q_1'}^{p_1'}(t)} <\infty$ for $t\in (0,\infty)$. By  \cite[Theorem 2.1]{YZ}, the PDE
\beq\label{NP} \pp_s u_s= - L_s  (u_s +f),\ \ s\in [0,t],\ u_t=0\end{equation}
  has a unique solution satisfying
  $$\|u\|_\infty+\|\nn u\|_\infty+\|\nn^2u\|_{\tt L_{q_1'}^{p_1'}(t)}+ \|(\pp_s +b^{(1)}  \cdot\nn) u\|_{\tt L_{q_1'}^{p_1'}(t)}<\infty,$$
  so that It\^o's formula (see \cite[Theorem 1.2.3]{RW24}) yields
  \beg{align*} \d \big\{u_r(X_{s,r}^{x})\big\}
  &= \big\{(L_r +\pp_r) u_r\big\}(X_{s,r}^{x})\d r
   +\big\<\si_r(X_{r}^{x})\d W_r, \nn u_r( X_{s,r}^{x})\big\>\\
   &= -L_r f (X_{s,r}^{x}) \d r +\big\<\si_r(X_{s,r}^{x})\d W_r, \nn u_r( X_{s,r}^{x})\big\>,\ \  r\in [s,t].\end{align*}
  Combining  this with $u_t=0, X_{s,s}^{x}=x$ and   \eqref{FW}, we derive
\beg{align*}&-u_s(x)= \E\big[u_t(X_{s,t}^{x})- u_s( X_{s,s}^{x})\big]
 =- \E\int_s^t L_r  f( X_{s,r}^{x})\d r\\
 &= - \int_s^t   P_{s,r} (L_r f)(x)\d r =- \int_s^t \pp_r P_{s,r} f(x)\d r
 = f(x)-P_{s,t} f(x).\end{align*}
 This together with \eqref{NP}   implies \eqref{BC} for $\ell'=1$.

Assume that \eqref{BC} holds for $\ell'=n$ for some $n\in \mathbb N$. Let  $u_s, \Theta_s$ be in the proof of Proposition \ref{L*1}, and let
\beq\label{CV} \bar P_{s,t}f(x):=  \E[f(\Theta_t(X_{s,t}^{\Theta_s^{-1}(x)})] = P_{s,t} (f\circ \Theta_t)(\Theta_s^{-1}(x)).\end{equation}
Since the coefficients $\bar{b},\bar{\sigma}$ in \eqref{XC} for the associated SDE to $\bar P_{s,t}$ satisfy $(C)$ for $\ell'=n$, we obtain
$$\pp_s \bar P_{s,t} f  = - \Big(\ff 12  {\rm tr} (\bar \si_s \bar \si_s^*\nn^2) +\bar b_s\cdot \nn \Big) \bar P_{s,t} f. $$
  This together with \eqref{PL0} and \eqref{CV} implies \eqref{BC}.

Now, by \eqref{BC} and   It\^o's formula,   we find a martingale $M_s$ such that
 \beg{align*} \d \big\{P_{s,t} f(\tt X_s^{x})\big\} &= [\big(\pp_s+ \tt L_s\big)
 P_{s,t} f](\tt X_s^{x})\d s+\d M_s\\
&= \big\<(\tt b_s^{0}-b_s^{0})(\tt X_s^{x}),\ \nn  P_{s,t} f(\tt X_s^{x})\big\>\d s+\d M_s,\ \ s\in [0,t],\end{align*}
which implies  \eqref{DH1'}.
\end{proof}






\end{document}